\documentclass[12pt,a4paper,leqno]{article}

\usepackage{amsfonts}
\usepackage{amssymb}
\usepackage{fullpage}
\usepackage{graphicx}

\newtheorem{num}{\indent\hskip-4.5pt}[section]
\renewcommand{\thenum}{\rm \arabic{section}.\arabic{num}}
\newcommand{\alku}{\begin{num}\hskip-6pt .\hskip5pt }
\newcommand{\loppu}{\end{num}}

\newcommand{\be}{\addtocounter{num}{1}\begin{equation}}
\newcommand{\ee}{\end{equation}}
\newcommand{\ots}[1]{\bigskip \refstepcounter{num} \indent
\thenum.\hskip 0.2truecm {\bf #1.} \hskip 0.4truecm}

\newcommand{\bea}{\addtocounter{num}{1}\begin{eqnarray}}  
\newcommand{\eea}{\end{eqnarray}}

\font\ff=eusm10 scaled 1200

\def\QM{\mathrm{QM}}
\def\Im{\mathrm{Im}\,}
\def\Re{\mathrm{Re}\,}
\def\M{\mathcal{M}}
\def\N{\mathbb{N}}
\def\Z{\mathbb{Z}}
\def\R{\mathbb{R}}
\def\C{\mathbb{C}}

\def\lem{{\bf Lemma.\hskip 0.5truecm}}
\def\pro{{\bf Proposition.\hskip 0.5truecm}}
\def\cor{{\bf Corollary.\hskip 0.5truecm}}
\def\theo{{\bf Theorem.\hskip 0.5truecm}}
\def\rem{{\bf Remark.\hskip 0.5truecm}}

\def\conj{{\bf Conjecture.\hskip 0.5truecm}}

\def\E{\hbox{\ff E}}
\def\K{\hbox{\ff K}}
\def\mut{\tilde{\mu}}
\def\phit{\tilde{\varphi}}

\def\proof{{\bf Proof.\hskip 0.5truecm}}       

\newcounter{minutes}
\divide\time by 60
\newcounter{hours}
\multiply\time by 60
\addtocounter{minutes}{-\time}

\begin{document}

\title{\bf GENERALIZED ELLIPTIC INTEGRALS}
\author{\bf Ville Heikkala\thanks{\fussy Author supported by the Magnus
Ehrnrooth fund of the Finnish Academy of Science and Letters.}, \\
\bf Mavina K.\ Vamanamurthy, \\
\bf and \\
\bf Matti Vuorinen}

\date{\begin{center}
\texttt{\small File:~\jobname .tex, 
        printed: \number\year-\number\month-\number\day, 
        \thehours.\ifnum\theminutes<10{0}\fi\theminutes}
\end{center}}

\renewcommand{\thefootnote}{\fnsymbol{footnote}}
\setcounter{footnote}{1}
\maketitle

\renewcommand{\thefootnote}{\arabic{footnote}}
\setcounter{footnote}{0}

\medskip
\begin{abstract}
Jacobi's elliptic integrals and elliptic functions arise
naturally from the Schwarz-Christoffel conformal transformation of
the upper half plane onto a rectangle. In this paper we study
generalized elliptic integrals which arise from the analogous mapping
of the upper half plane onto a quadrilateral and obtain sharp 
monotonicity and convexity properties for certain combinations
of these integrals, thus generalizing analogous well-known results
for classical conformal capacity and quasiconformal distortion 
functions. An algorithm for the computation of the modulus of the
quadrilateral is given.
\end{abstract} 

\bigskip

{\bf 2000 Mathematics Subject Classification:} 
Primary 33B15, 33C05, Secondary 30C62.

\bigskip


\section{Introduction}

Given complex numbers
$a,b,$ and $c$ with $c\neq 0,-1,-2, \dots $,
the {\em Gaussian hypergeometric function} is the analytic
continuation to the slit plane $\C \setminus [1,\infty)$ of
the series
\be \label{eq:hypdef}
F(a,b;c;z) = {}_2 F_1(a,b;c;z) = 
\sum_{n=0}^{\infty} \frac{(a,n)(b,n)}{(c,n)} \frac{z^n}{n!}\,,\:\: |z|<1 \,.
\ee
Here $(a,0)=1$ for $a \neq 0$, and $(a,n)$
is the {\em shifted factorial function}
or the {\em Appell symbol}
$$
(a,n) = a(a+1)(a+2) \cdots (a+n-1)
$$
for $n \in \N \setminus \{0\}$, where 
$\N = \{ 0,1,2,\ldots\}$.

A {\em generalized modular equation of order (or degree)} $p>0$ is
\be \label{eq:abc}
\frac{F(a,b;c;1-s^2)}{F(a,b;c;s^2)} = p
\frac{F(a,b;c;1-r^2)}{F(a,b;c;r^2)}\,, \:\: 0 < r < 1\,.
\ee
Sometimes we just call this an $(a,b,c)$-modular equation of order $p$
and we usually
assume that $a,b,c>0$ with $a+b \ge c$, in which case this equation uniquely
defines $s$, see Lemma \ref{muabcmonot}. 

Many particular cases of (\ref{eq:abc}) have been studied in the
literature on both
analytic number theory and geometric function theory, \cite{Be}, 
\cite{BB}, \cite{BBG}, \cite{AVV}, \cite{AQVV}, \cite{LV}. 
The classical case $(a,b,c)= (\frac{1}{2}, \frac{1}{2}, 1)$
was studied already by Jacobi and many others in the nineteenth century, see
\cite{Be}.     
In 1995 B.\ Berndt, S.\ Bhargava, and F.\ Garvan published
an important paper \cite{BBG} in which they studied the case 
$(a,b,c)= (a,1-a,1)$ and $p$ an integer. For several rational
values of $a$ such as $a=\frac{1}{3}, \frac{1}{4}, \frac{1}{6}$ 
and integers $p$ (e.g.\ $p= 2,3,5,7,11,...$)  
they were able to give proofs for numerous algebraic
identities stated by Ramanujan in his unpublished notebooks. 
These identities involve $r$ and $s$ from (\ref{eq:abc}).
After the publication of \cite{BBG} many papers have appeared on
modular equations, see e.g.\ \cite{AQVV}, \cite{Be},
\cite{CLT}, \cite{Q}, and \cite{S}. 

To rewrite (\ref{eq:abc}) in a slightly shorter form, we use the decreasing
homeomorphism $\mu_{a,b,c} : (0,1) \to (0, \infty)$, defined by
\be \label{eq:muacdef}
\mu_{a,b,c}(r) = \frac{B(a,b)}{2} 
\frac{F(a,b;c;{r'}^2)}{F(a,b;c;r^2)}\,, \:\: r \in (0,1)
\ee
for $a, b, c >0$, $a+b \ge c$, where $B$ is 
the beta function, see (\ref{eq:betadef}) below.  
We call $\mu_{a,b,c}$ the {\em generalized modulus}, cf.\ \cite[(2.2)]{LV}.
We can now write (\ref{eq:abc}) as
\be \label{newmodeq}
\mu_{a,b,c}(s) = p \thinspace \mu_{a,b,c}(r)\,, \:\: 0<r<1\,.
\ee

With $p = 1/K$, $K > 0$, the solution of (\ref{eq:abc}) is then given by
\be \label{eq:phiackdef}
s= \varphi_{K}^{a,b,c}(r) = \mu_{a,b,c}^{-1}(\mu_{a,b,c}(r)/K)\,.
\ee
We call $\varphi_K^{a,b,c}$ 
{\em the $(a,b,c)$-modular function with degree $p= 1/K$}
\cite{BBG}, \cite[(1.5)]{AQVV}.

In the case $a<c$ we also use the notation
$$
\mu_{a,c} = \mu_{a,c-a,c}\,, \:\: 
\varphi_{K}^{a,c} = \varphi_{K}^{a,c-a,c}\,.
$$

For $0 < a < \min \{c,1\}$ and $0 < b < c \le a+b $, define the 
{\it generalized complete elliptic integrals 
of the first and second kinds} 
(cf.\ \cite[(1.9), (1.10), (1.3), and (1.5)]{AQVV}) on $[0,1]$ by
\be \label{eq:Kdef}
\K = \K_{a,b,c} = \K_{a,b,c}(r) = \frac{B(a,b)}{2} F(a,b;c;r^2)\,,
\ee
\be \label{eq:Edef}
\E = \E_{a,b,c} = \E_{a,b,c}(r) = \frac{B(a,b)}{2} F(a-1,b;c;r^2)\,,
\ee
\be \label{eq:KpEpdef}
\K' = \K_{a,b,c}' = \K_{a,b,c}(r')\,, \:\: 
\mbox{\rm and} \:\: \E' = \E_{a,b,c}' = \E_{a,b,c}(r')
\ee
for $r \in (0,1)$, $r' = \sqrt{1-r^2}$. The end values are defined by limits
as $r$ tends to $0+$ and $1-$, respectively.
In particular, we denote $\K_{a,c} = \K_{a,c-a,c}$ and 
$\E_{a,c} = \E_{a,c-a,c} \,.$ 
Thus, by (\ref{eq:hypas}) below,
$$
\K_{a,b,c}(0) = \E_{a,b,c}(0) = \frac{B(a,b)}{2}
$$
and
$$
\E_{a,b,c}(1) = \frac{1}{2}
\frac{B(a,b) B(c,c+1-a-b)}{B(c+1-a,c-b)}\,,
\:\: \K_{a,b,c}(1) = \infty\,.
$$
Note that the restrictions on $a,b$ and $c$ ensure 
that the function $\K_{a,b,c}$ is increasing and unbounded
whereas $\E_{a,b,c}$ is decreasing and bounded,
as in the classical case $a = b = \frac{1}{2}, c=1$.
Note also that our terminology differs from that of \cite[Section 5.5]{BB}, 
where generalized elliptic integrals refer to the particular case
$c=1$. 

In this paper we study the modular function $\varphi_K^{a,b,c}$ and the
generalized modulus $\mu_{a,b,c}$ as
well as the generalized elliptic integrals $\K_{a,b,c}$ and $\E_{a,b,c}$. 
In the case $b = 1-a$, $c = 1$, these functions 
coincide with the special cases $\varphi_K^a$, $\mu_a$, $\K_a$, and
$\E_a$ which were studied in \cite{AQVV}.

In Section \ref{sect:mapping} we construct a conformal mapping
from a quadrilateral with internal 
angles $b \pi$, $(c-b) \pi$, $(1-a) \pi$, and $(1-c+a) \pi$ onto the
upper half plane. We denote this mapping by $\mbox{\rm sn}_{a,b,c}$.
If $b = 1-a$ and $c = 1$, 
this mapping reduces to the generalized elliptic
sine function $\mbox{\rm sn}_{a}$ in \cite[2.1]{AQVV}.
We given an algorithm for the computation of the conformal modulus
of a quadrilateral and its implementation in the 
$\mbox{\rm Mathematica}^{\circledR}$ language.
In a later work, we study the dependence of the modulus on the
geometry of the quadrilateral, \cite{DV}, \cite{RV}.

In Section \ref{sect:special} we recall some basic properties of the
hypergeometric, gamma, and beta functions, that are used in the sequel. 

Section \ref{sect:modular} contains our main results:
differentiation formulas and
monotone properties of the generalized
elliptic integrals and of several combinations of these functions.
Theorems \ref{monotquot} and \ref{monpolquot}, which provide 
sufficient conditions for the
quotients of two Maclaurin series to be monotone increasing, seem to be of
independent interest and, as far as we know, these results are new.
Also a similar result is given for the quotient of two polynomials
of the same degree.

The special functions $\K_{a,b,c}, \mu_{a,b,c}$ and quotients of
hypergeometric functions are presently studied
very intensively. The manuscript \cite{HLVV} is a direct continuation
of this work and for some particular parameter triples $(a,b,c)$ there are
very recent results by many authors \cite{B1}-\cite{B4}, \cite{KS},
\cite{WZQC}, \cite{ZWC}.

Throughout this paper we denote $r' = \sqrt{1-r^2}$ whenever
$r \in (0,1)$. The standard symbols
$\C, \R, \Z$, and $\N$ denote the sets of complex numbers, real numbers, 
integers, and natural numbers (with zero included), respectively.

\bigskip


\section{The Schwarz-Christoffel map onto a quadrilateral} 
\label{sect:mapping}

For  $0 < a,b < 1, \max\{ a+b, 1 \}  \le c \le \min \{a,b\} + 1$,
$r \in (0,1)$, let
$g_r(t) = t^{b-1}(1-t)^{c-b-1}(1-r^2t)^{-a}$, $t \in \C$, $\Im t \ge 0$, 
denote the analytic branch
for which the argument of each of the factors $t$, $1-t$, and 
$1-r^2 t$ is $\pi$ whenever it is real and negative. 
Denote $C = C(b,c) = 1/B(b,c-b)$.
We define the generalized Jacobi
sine function $\mbox{\rm sn}_{a,b,c}(w) = \mbox{\rm sn}_{a,b,c}(w,r)$
as the inverse of the function $f$
given on the closed upper half plane by
\bea
w = f(z) = f_{a,b,c}(z) & = & C \int_{0}^{z} g_r(t) dt \nonumber \\ 
 & = & C \int_{0}^{z} 
t^{b-1}(1-t)^{c-b-1}(1-r^2t)^{-a} dt \label{eq:ellsin} \\ 
 & = &  e^{i(a+b+1-c)\pi} C r^{-2a} 
\int_{0}^{z} t^{b-1}(t-1)^{c-b-1}(t-1/r^2)^{(1-a)-1} dt\,. \nonumber
\eea
Recall the Euler integral representation \cite[Theorem 2.2.1]{AAR}  
\cite[15.3.1]{AS}
\be \label{eq:Fint}
F(a,b;c;z) = 
\frac{\Gamma(c)}{\Gamma(b)\Gamma(c-b)} 
\int_{0}^{1} t^{b-1}(1-t)^{c-b-1}(1-tz)^{-a} dt
\ee
$$
= C(b,c) e^{i(a+b+1-c)\pi} \int_0^1 t^{b-1}(t-1)^{c-b-1}(tz-1)^{-a} dt
$$
for $\Re c > \Re b > 0$ and
$z \in \C \setminus \{u \in \R \: | \: u \ge 1 \}$. 

The next result is well-known, but for the sake of completeness a proof
is given. As general references concerning
 the Schwarz-Christoffel mapping \cite{M} and \cite{N} may be mentioned.

\bigskip
\alku \pro \label{nehariconv}   
The Schwarz-Christoffel mapping of the upper half plane $H$ onto 
a quadrilateral
with turning angles $\beta_k \pi$ at the vertices $a_k \in \R$ is
$$ w= F(z)= \int_0^z  \frac{dz}{\Pi_{k=1}^4(z-a_k)^{\beta_k}},\quad
-1<\beta_k<1, \sum_{k=1}^4 \beta_k = 2, $$
where $a_k$ are the points on the real axis that $F$ maps to the four
vertices of the quadrilateral.
\loppu

\proof As in \cite[pp.\ 192-193]{N} let 
$T:H \to U$ be the M\"obius transformation $ T(z) = (z-i)/(z+i)\,,$ where  $U$
is the unit disk $ \{ z:  |z|<1\}. $
Then $f(z) = F(T^{-1}(z))$ maps 
$U$ onto the quadrilateral $Q \,$ and
$$ f(\zeta) = \int_0^{\zeta} \frac{dw}{\Pi_{k=1}^4(w-\zeta_k)^{\beta_k}}\,,$$
where $\zeta_k$ are the points on the unit circle that map onto
the four vertices of the quadrilateral $Q$. Now
$$
\frac{\zeta f''(\zeta)}{f'(\zeta)} = 
-\sum_{k=1}^4  \beta_k \frac{\zeta}{(\zeta-\zeta_k)} = 
\sum_{k=1}^4  \beta_k \frac{1+ \overline{\zeta_k {} } \zeta }
{1- \overline{\zeta_k {} } \zeta }
$$
and recalling that the $\beta_k$ sum to $2$ we have that
$$
1+\frac{\zeta f''(\zeta)} {f'(\zeta)} =
\frac{1}{2} \sum_{k=1}^4 \beta_k 
\frac{1+ \overline{\zeta_k {} } \zeta } {1- \overline{\zeta_k {} } \zeta } \,.
$$
For all points $\zeta\in U$
$$
\Re  \left\{  1 +\frac{\zeta f''(\zeta)}{f'(\zeta)}  \right\} = 
\frac{1}{2} \sum_{k=1}^4 \beta_k   
\Re \left\{ \frac{1+  \overline{\zeta_k {} } \zeta  }
{1-  \overline{\zeta_k {} } \zeta } \right\}= 
\frac{1}{2} \sum_{k=1}^4 \beta_k    
\frac{1-| \zeta|^2  }{|1-  \overline{\zeta_k {} } \zeta |^2} >0 \,.
$$
This proves that $f$ maps  $U$ bijectively onto a convex domain, since
$\Re \{1+ {\zeta f''(\zeta)}/{f'(\zeta)}  \}>0$ is a
necessary and sufficient condition for this to be true, 
see \cite[Ex.\ 5, p.\ 224]{N}. $\qquad \square$  

\bigskip

\begin{figure}[htbp]
\centering
\includegraphics[width=70mm]{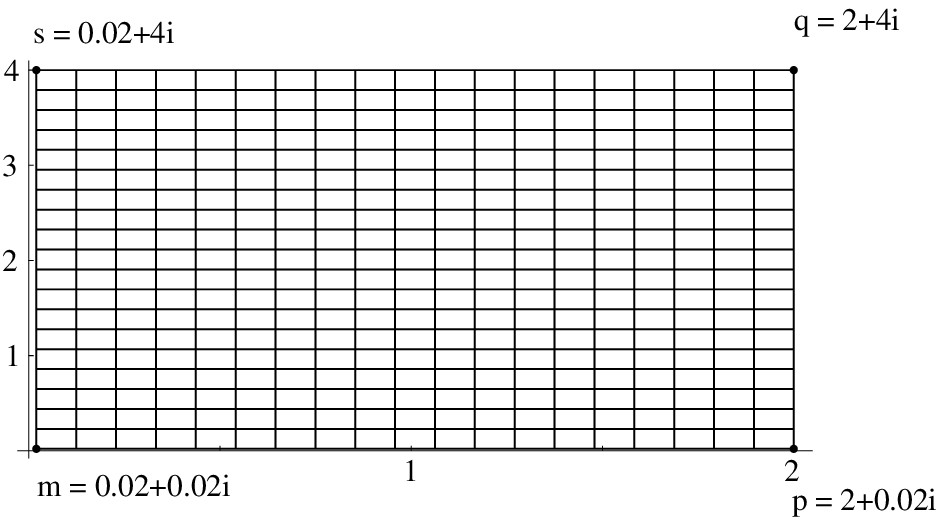}
\hfill
\includegraphics[width=85mm]{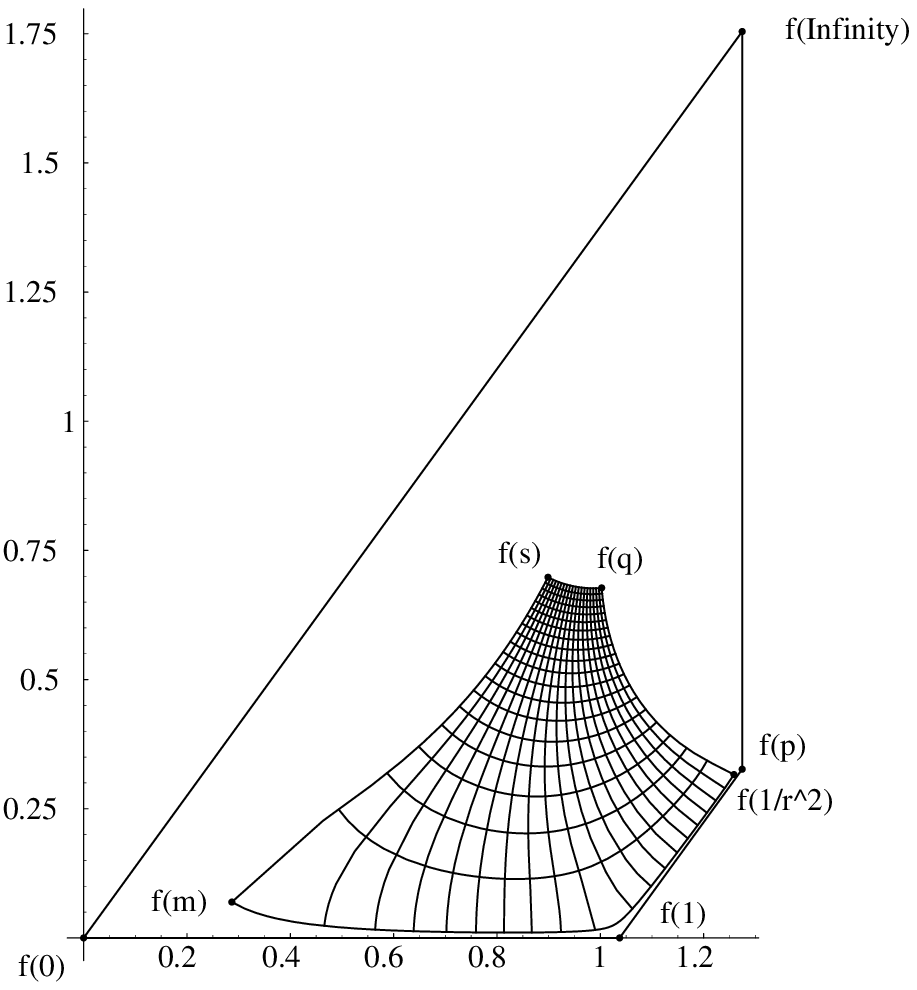}
\caption{The image quadrilateral and the image of a grid under the mapping
$f_{a,b,c}$ with $a = 0.2, b = 0.3, c = 1.0, r = 0.7$}
\end{figure}

\bigskip

\alku \theo \label{ellsintheo}
Let $H$ denote the closed upper half-plane $\{z \in \C \: | \: 
\Im z \ge 0\}$ and let  
$0 < a, b < 1, \max \{a+b,1 \} \le c \le 1 + \min \{a,b \}$,
$r \in (0,1)$. Then the
function $f$ in (\ref{eq:ellsin}) is a homeomorphism of $H$ onto the
quadrilateral $Q$ with vertices 
$$
f(0) = 0\,,\:\: f(1) = F(a,b;c;r^2)\,,
$$
$$
f(1/r^2) = f(1) + 
\frac{B(c-b,1-a)}{B(b,c-b)} e^{(b+1-c) i \pi} {r'}^{2(c-a-b)}  
F(c-a,c-b;c+1-a-b;{r'}^2)\,,
$$
and
$$
f(\infty) = f(1/r^2) 
+ \frac{B(1-a,a+1-c)}{B(b,c-b)}  e^{i(a+b+1-c)\pi} r^{2(1-c)} {r'}^{2(c-a-b)} 
F(1-b,1-a;2-c;r^2)\,,
$$
and interior angles $b\pi , (c-b)\pi, (1-a)\pi$ and $(a+1-c)\pi$, 
respectively, at these vertices. It is conformal in the interior of $H$.
\loppu

\medskip

\proof It is clear that $f(0) = 0$ and by (\ref{eq:Fint}) 
$f(1) = F(a,b;c;r^2)$. Next we evaluate
\begin{eqnarray*}
f(1/r^2) & = & C \int_{0}^{1/r^2} g_r(t) dt \\
 & = & C \int_{0}^{1} g_r(t) dt + C \int_{1}^{1/r^2} g_r(t) dt \\
 & = & F(a,b;c;r^2) + C \int_{1}^{1/r^2} g_r(t) dt\,.
\end{eqnarray*}
To evaluate the second integral above, we make the change of variable
$t = 1/(1-{r'}^2 u)$ for which $dt = (1-{r'}^2 u)^{-2} {r'}^2 du$.
We observe that this change of variable is simply the restriction to
reals of the plane M\"obius transformation taking the ordered triple
$(1, 1/r^2, 0)$ onto the ordered triple $(0, 1, \infty)$.
Then
\begin{eqnarray*}
g_r(t) dt & = & (1-{r'}^2 u)^{1-b} 
\left(1-\frac{1}{1-{r'}^2 u} \right)^{c-b-1}
\left(1-\frac{r^2}{1-{r'}^2 u} \right)^{-a}
(1-{r'}^2 u)^{-2} {r'}^2 du \\
 & = & (-1)^{b+1-c} {r'}^{2(c-a-b)} u^{c-b-1} (1-u)^{-a}
(1-{r'}^2 u)^{a-c} du \\
 & = & e^{(b+1-c) i \pi} {r'}^{2(c-a-b)} u^{c-b-1} (1-u)^{-a}
(1-{r'}^2 u)^{a-c} du
\end{eqnarray*}
and by (\ref{eq:Fint}) we get
\begin{eqnarray*}
C \int_{1}^{1/r^2} g_r(t) dt & = & 
e^{(b+1-c) i \pi} {r'}^{2(c-a-b)} 
\frac{\Gamma(c)}{\Gamma(b)\Gamma(c-b)} 
\int_{0}^{1} u^{c-b-1} (1-u)^{-a}(1-{r'}^2 u)^{a-c} du \\
 & = & e^{(b+1-c) i \pi} {r'}^{2(c-a-b)} 
\frac{\Gamma(c)}{\Gamma(b)\Gamma(c-b)} 
\frac{\Gamma(c-b)\Gamma(1-a)}{\Gamma(c+1-a-b)} \cdot \\
& & \cdot F(c-a,c-b;c+1-a-b;{r'}^2) \\
 & = & \frac{B(c-b,1-a)}{B(b,c-b)} 
 e^{(b+1-c) i \pi} {r'}^{2(c-a-b)} 
 F(c-a,c-b;c+1-a-b;{r'}^2) \\
 & = & \frac{B(c-b,1-a)}{B(b,c-b)} 
 e^{(b+1-c) i \pi} r^{2(1-c)} {r'}^{2(c-a-b)} \cdot \\ 
 & & \cdot F(1-a,1-b;c+1-a-b;{r'}^2)\,,
\end{eqnarray*}
where the last expression follows from (\ref{eq:hypas}) below.
Hence $f(1/r^2)$ has the value claimed.
We proceed to evaluate 
the remaining value, namely
$$
f(\infty) = C \int_{0}^{\infty} g_r(t) dt
= C \int_{0}^{1/r^2} g_r(t) dt + C \int_{1/r^2}^{\infty} g_r(t) dt\,.
$$
The first integral above equals $f(1/r^2)$.
To compute the second one, we apply the change of variable
$t = (1-r^2 v)/(r^2 (1-v))$. We observe that
this change of variable is simply the restriction to reals of the
plane M\"obius transformation taking the ordered triple 
$(1/r^2, \infty, 1)$ onto the ordered triple $(0, 1, \infty)$. Then
$dt = (1-r^2)/(r^2(1-v)^2) dv$, $t = 1/r^2$ gives $v = 0$, and
$t = \infty$ gives $v = 1$.
We get
\begin{eqnarray*}
g_r(t) dt & = & \left(\frac{1-r^2 v}{r^2 (1-v)} \right)^{b-1}
\left(1-\frac{1-r^2 v}{r^2 (1-v)} \right)^{c-b-1} \cdot \\
 & & \cdot
\left(1-\frac{r^2 (1-r^2 v)}{r^2 (1-v)} \right)^{-a}
\left(\frac{1-r^2}{r^2 (1-v)^2}\right) dv \\
 & = & r^{2(1+a-b+b+1)} (1-v)^{a-c} (1-r^2 v)^{b-1}
(r^2 - 1)^{c-b-1} (r^2 (r^2-1)v)^{-a} (1-r^2) dv \\
 & = & (-1)^{a+b+1-c} {r'}^{2(c-a-b)} r^{2(1-c)} v^{-a} (1-v)^{a-c}
 (1-r^2v)^{b -1} dv
\end{eqnarray*}
so that
\begin{eqnarray*}
C \int_{1/r^2}^{\infty} g_r(t) dt & = & 
(-1)^{a+b+1-c} {r'}^{2(c-a-b)}  r^{2(1-c)} C 
\int_{0}^{1} v^{-a} (1-v)^{a-c} (1-r^2 v)^{b-1} dv \\
 & = & (-1)^{a+b+1-c} {r'}^{2(c-a-b)}  r^{2(1-c)} 
\frac{B(1-a, 1+a-c)}{B(b,c-b)} \cdot \\
 & & \cdot F(1-b, 1-a; 2-c; r^2)\,.
\end{eqnarray*}
The claimed value for $f(\infty)$ follows.

It follows from the
formula (\ref{eq:ellsin}) that (see e.g.\ \cite[pp.\ 128--134]{M}) $f$ is a 
Schwarz-Christoffel transformation which maps $H$ onto a quadrilateral $Q$ with
vertices $f(0)$, $f(1)$, $f(1/r^2)$, and $f(\infty)$ and interior angles
$b \pi$, $(c-b) \pi$, $(1-a) \pi$, and $(1-c+a) \pi$ in counterclockwise
order. 
By Proposition \ref{nehariconv} $f$ is univalent.
$\qquad \square$

\bigskip

\alku \cor \label{modquad}
Let $0 < a, b < 1, \max \{a+b,1 \} \le c \le 1 + \min \{a,b\}$, 
and let $Q$ be 
a quadrilateral in the upper half plane with vertices $0, 1, A$ and $B$, 
the interior angles at which are, respectively, $b\pi, (c-b)\pi, 
(1-a)\pi$ and $(1+a-c)\pi$. Then the conformal modulus of $Q$ 
(cf.\ \cite{LV}) is given by
$$
\mbox{\rm mod}(Q) = \K(r') / \K(r),
$$ 
where $r \in (0,1)$ satisfies the equation
\be \label{eq:star}
A - 1 = \frac{L {r'}^{2(c-a-b)} F(c-a,c-b;c+1-a-b; {r'}^2)}{F(a,b;c;r^2)} =
G(r) \, ,
\ee
say, and
$$
L = \frac{B(c-b,1-a)}{B(b,c-b)} e^{(b+1-c)i \pi}.
$$

\loppu

\medskip

\proof
Clearly, $\arg(A-1) = (b+1-c)\pi = \arg(L)$. Since 
$G(0) = \infty$ and $G(1) = 0$, 
it follows that a unique solution $r \in (0,1)$ of equation 
(\ref{eq:star}) exists. 
Let $f$ be as in Theorem \ref{ellsintheo} and let $g = f / f(1)$. 
Then $g$ maps the upper half plane $H$ onto $Q$, with
$g(0) = 0$, $g(1) = 1$, $g(1/r^2) = A$, and $g(\infty) = B$. The function 
$h = \mbox{\rm sn}^{-1}$ maps the first quadrant conformally \cite{Bo} 
onto the standard rectangle $R$, with 
$h(0)=0$, $h(1) = \K(r)$, $h(1/r) = \K(r) + i \K(r')$, and 
$h(\infty) = i \K(r')$. 
Hence the function $k = h(\sqrt{})$ maps $H$ conformally onto $R$. Thus, by 
conformal invariance,
$mod(Q) = \K(r') /\K(r)$. $\qquad\square$

\bigskip

\alku \rem \label{quadrspecc}
{\rm 
(1) The quadrilateral $Q$ in Theorem \ref{ellsintheo} 
reduces to a trapezoid if and only if $c=1$ or $c=a+b$,
to a parallelogram if and only if $c=1=a+b$, \cite{AQVV} and 
to a rectangle (the classical case)
if and only if $a=b=\frac{1}{2}$ and $c=1$, \cite{Bo}.

(2) The hypotheses in Corollary \ref{modquad} imposed on the triple
$a,b,c$ imply that the quadrilateral $Q$ is convex.
}
\loppu

\medskip

\alku \rem \label{bowmanmod}
{\rm 
Bowman \cite[pp.\ 103-104]{Bo} gives a formula for the conformal
modulus of the quadrilateral with vertices $0$, $1$, $1+hi$, and $(h-1)i$ 
when $h>1$ as $q=\K(r)/\K(r')$ where
$$ 
r = \left(\frac{t_1-t_2}{t_1+t_2}\right)^2 \, , \quad 
t_1 = \mu^{-1}\left(\frac{\pi}{2c}\right) \, , \quad
t_2 = \mu^{-1}\left(\frac{\pi c}{2}\right) \, , \quad
c = 2h - 1 \, .
$$
Therefore, the quadrilateral can be conformally mapped onto the
rectangle $0$, $1$, $1+qi$, $qi$ with vertices corresponding.
}
\loppu

\bigskip

\alku {\bf Computational discovery.} \label{compdiscov}
{\rm
We have written a 
$\mbox{\rm Mathematica}^{\circledR}$ function that computes the modulus
of the quadrilateral with vertices at $0, 1, A, B$ where 
$\Im A >0, \Im B > 0$. This led to the following discovery
for symmetric quadrilaterals:
\be \label{eq:symmquad}
\mathrm{If} \:\: |B|=1 \:\: \mathrm{and} \:\: 
2{\mathrm{arg}} A = {\mathrm{arg}} B, \:\: 
\mbox{\rm then the modulus is equal to} \:\: 1 \,.
\ee
It is not difficult to prove this analytically 
(see e.g.\ \cite[p.\ 433]{Hen} or \cite{Her}).
The following
figure illuminates the variation of the modulus of the
quadrilateral with vertices $0, 1, x+i y, i$ in the first
quadrant.
}
\loppu

\bigskip

\begin{figure}[htbp]
\centering
\includegraphics[width=100mm]{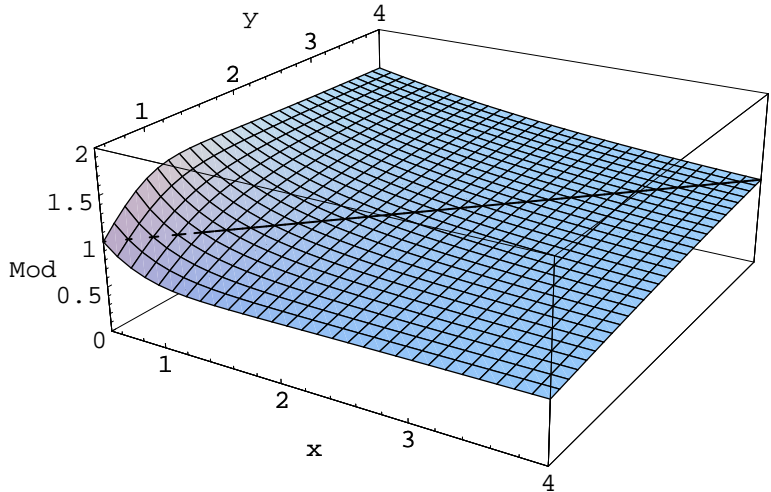}
\caption{Modulus of the quadrilateral with vertices $0, 1, x+iy, i$ and
the line $(x,x,1)$.}
\end{figure}

\bigskip

\alku {\bf Duplication formula for quadrilateral modulus.} \label{quadrdup}
{\rm
For the sake of convenient notation, denote by 
$\QM(a,b,c,d)$ the modulus of the polygon with vertices $a,b,c,d$. 
That is, $\QM(a,b,c,d) = q$ if and only if there exists a conformal map 
of the polygon onto the rectangle with vertices $1+iq, iq, 0, 1$ such 
that the vertices correspond. Clearly $\QM(1+iq, iq, 0, 1) = q$. 
Remark \ref{bowmanmod}
gives now an explicit formula for $\QM(1+ih, i(h-1), 0, 1)$ when $h > 1$.
(By the way, it seems likely that $\QM(1+ih, i(h-1), 0, 1) \in [h-1,h]$.)
By (\ref{eq:symmquad}) we see that $\QM(te^{i\pi/4},i,0,1) = 1$ 
for all $t > 0$ and it is clear by symmetry that
$$
\QM(te^{i(\pi/4-\alpha)},i,0,1) = 1/\QM(te^{i(\pi/4+\alpha)},i,0,1) 
$$
for all $\alpha \in (0, \pi/4)$ and all $t > 0$. From 
\cite[Corollary 2.3]{AQVV} we have an expression for the parallelogram case
$\QM(1+r e^{i\alpha}, r e^{i\alpha}, 0, 1)$.

We are not familiar with cases other than those mentioned above
where the values of $\QM(a,b,c,d)$ could be expressed in reasonably
simple form. Therefore it might be of some interest to record a
duplication formula for $\QM(a,b,c,d)$ which follows from an elementary
symmetry consideration. To this end, let $h,k>0$ and consider
the quadrilaterals $(1+ih,ik,-ik,1-ih)$ and $(1+ih,ik,0,1)$. 
It is clear by symmetry that
$$
\QM(1+ih,ik,-ik,1-ih) = 2 \QM(1+ih,ik,0,1) \,.
$$
We now use the invariance of $\QM(a,b,c,d)$ under homotheties to
express this result in a simpler form. Choose a homothety $T(z)=az+b$
with $T(-ik)=0, T(1-ih)=1$. Then $a=1/c, b=ik/c$ with $c=1+i(k-h)$.
Because
$$ 
T(1+ih)= (1+i(h+k))/c\,, \: T(ik)= 2ik/c\,, \:
T(-ik)=0\,, \: \mathrm{and} \: T(1-ih)=1\,,  
$$
we have for $h,k>0, c=1+i(k-h)$ that
\be \label{eq:qmdupl}
\QM((1+i(h+k))/c, 2ik/c, 0, 1)= 2\QM(1+ih, ik, 0, 1)\,. 
\ee
}
\loppu

\bigskip

\begin{table}[htbp]
\centering
\begin{tabular}{|c|c|c|c|c|c|}
\hline
m $\backslash$ n & 1 & 2 & 3 & 4 & 5 \\
\hline
1 & 1.000000 & 1.279261 & 1.354244 & 1.383086 & 1.397799 \\
2 & 0.781700 & 1.000000 & 1.127663 & 1.201627 & 1.248066 \\
3 & 0.738419 & 0.886789 & 1.000000 & 1.080783 & 1.138566 \\
4 & 0.723020 & 0.832204 & 0.925254 & 1.000000 & 1.058739 \\
5 & 0.715410 & 0.801239 & 0.878297 & 0.944519 & 1.000000 \\
\hline
\end{tabular}
\caption{Numerical values of $\QM(m+in,i,0,1)$ for 
$m,n = 1,\ldots,5$. The values have been truncated 
to six decimal places} \label{tab:QMvals}
\end{table}

Below we give a $\mbox{\rm Mathematica}^{\circledR}$ program 
for the computation
of $QM(A,B,0,1).$ This program computes in fact the values
given in the attached table.

\begin{verbatim}
(*FILE:  qm4hvv.m begins 2007-03-18  *)
(* USAGE:  <<qm4hvv.m  *)
H2F1[a_,b_,c_,z_] := Hypergeometric2F1[a,b,c,z];
KK[r_]:= (Pi/2)Hypergeometric2F1[1/2,1/2,1,r*r];
myTmp[a_,b_,c_,x_]:=((1-x*x)^(c-a-b))*
     H2F1[c-a,c-b,c+1-a-b,1-x^2]/H2F1[a,b,c,x^2];
QM[A_,B_]:=
Module[{b=Arg[B]/Pi,c=(Pi-Arg[A-1]+Arg[B])/Pi//Simplify,
   a=1-(Arg[A-1]-Arg[A-B])/Pi//Simplify,L,r,t,f,x,init},
   L=Beta[c-b,1-a]/Beta[b,c-b]*Exp[I*(b+1-c)*Pi];
   f[x_]:=myTmp[a,b,c,x]-Re[(A-1)/L];
   r=x/.FindRoot[f[x]==0,{x,10^-6},WorkingPrecision->50, 
     MaxIterations->50];
	KK[Sqrt[1-r*r]]/KK[r]];
Do[Do[Print["QM[",m,"+I*",n,",I]=",QM[m+I*n,I]], {m,1,5}], {n,1,5}]
(*FILE:  qm4hvv.m ends *)
\end{verbatim}

This program enables one to carry out experiments with
the moduli of quadrilaterals.
The dependence of the modulus of a polygonal quadrilateral 
 on the geometry is the subject of our subsequent work
  \cite{DV}, \cite{RV}. 

\bigskip


\section{Hypergeometric functions}
\label{sect:special}

Let $\Gamma$ denote Euler's {\em gamma function} 
and let $\Psi$ be its logarithmic derivative (also called the
{\em digamma function}),
$\Psi(z) = \Gamma'(z)/\Gamma(z)$. By \cite[p.\ 198]{Ah} the function
$\Psi$ and its derivative have the series expansions
\be \label{eq:psiser}
\Psi(z) = -\gamma - \frac{1}{z} + \sum_{n=1}^{\infty} \frac{z}{n(n+z)}\,,
\:\: 
\Psi'(z) = \sum_{n=0}^{\infty} \frac{1}{(n+z)^2}\,,
\ee
where $\gamma = -\Psi(1)= \lim_{n \to \infty}(\sum_{k=1}^{n}1/k - \log n) =
0.57721\ldots$ is the {\em Euler-Mascheroni constant}. From
(\ref{eq:psiser}) it is seen that $\Psi$ is strictly increasing on 
$(0, \infty)$ and that $\Psi'$ is strictly decreasing there, so that
$\Psi$ is concave. Moreover,
$\Psi(z+1)= \Psi(z)+ 1/z$ and $\Psi(\frac{1}{2})= -\gamma -2 \log 2$, 
see \cite[Ch.\ 6]{AS}.

For all $z \in \C \setminus \{0, -1, -2, \ldots\}$ 
and for all $n \in \N$ we have 
\be \label{eq:gamapp}
\Gamma(z+n) = (z,n) \Gamma(z)\,,
\ee
a fact which follows by induction \cite[12.12]{WW}. 
This enables us to extend the Appell symbol for all complex values 
of $a$ and $a+t$, except for non-positive integer values, by
\be \label{eq:appext}
(a,t) = \frac{\Gamma(a+t)}{\Gamma(a)}\,.
\ee
Furthermore, the gamma function satisfies the
reflection formula \cite[12.14]{WW}
\be \label{eq:gamref}
\Gamma(z)\Gamma(1-z) = \frac{\pi}{\sin(\pi z)}
\ee
for all $z \not\in \Z$. In particular, 
$\Gamma(\frac{1}{2})= \sqrt{\pi}$.

The {\em beta function} is defined for $\Re x > 0, \: \Re y > 0$ by
\be \label{eq:betadef}
B(x,y) = \int_{0}^{1} t^{x-1} (1-t)^{y-1} dt
= \frac{\Gamma(x) \Gamma(y)}{\Gamma(x+y)}\,.
\ee

We use the standard notation for contiguous hypergeometric functions
(cf.\ \cite{R1})
$$
F = F(a,b;c;z), \:\: F(a+) = F(a+1,b;c;z), \:\: F(a-) = F(a-1,b;c;z)\,,
$$
etc. We also let 
$$
v = v(z) = F\,, \:\: u = u(z) = F(a-)\,, \:\: v_1 = v_1(z) = v(1-z)\,,
\:\:\mbox{\rm and} \:\: u_1 = u_1(z) = u(1-z)\,.
$$
The derivative of $F$ can be written in the following several 
different forms which will be useful in deriving the fifteen important
Gauss contiguous relations \cite{R1}
\bea
\frac{dv}{dz} = \frac{dF}{dz} & = & 
\frac{a}{z} (F(a+)-F) =
\frac{b}{z} (F(b+)-F) = \frac{c-1}{z}(F(c-)-F) \nonumber \\
& = & 
\frac{ab}{c}F(a+, b+; c+) = 
\frac{1}{1-z}\left((a+b-c)F + \frac{(c-a)(c-b)}{c}F(c+)\right) 
\label{eq:hypder} \\
& = & \frac{(c-a)u+(a-c+bz)v}{z(1-z)} \nonumber
\eea
and
\bea
\frac{du}{dz} = \frac{dF(a-)}{dz} & = &
(a-1)\left(F+\frac{b-c}{c}F(c+)\right) =
 \frac{a-1}{z}(v-u)\,. \label{eq:hypder2}
\eea
In particular, from (\ref{eq:hypder}) it follows that
(cf.\ \cite[Theorem 3.12 (3)]{AQVV})
\be \label{eq:fuv}
\frac{ab}{c} z(1-z) F(a+1, b+1; c+1; z) = 
(c-a) u(z) + (a-c+bz)v(z)\,.
\ee


The behavior of the
hypergeometric function near $z = 1$ in the three cases 
$\Re(a+b-c) < 0$, $a+b = c$, and $\Re(a+b-c) > 0$, 
respectively, is given by
\be \label{eq:hypas}
\left\{ \begin{array}{l}
F(a,b;c;1) = 
\frac{\Gamma(c) \Gamma(c-a-b)}{\Gamma(c-a) \Gamma(c-b)}\,, \\[.5cm]
B(a,b)F(a,b;a+b;z)+\log(1-z) = 
R(a,b)+{\rm O}((1-z)\log(1-z))\,, \\[.4cm]
F(a,b;c;z) = (1-z)^{c-a-b} F(c-a,c-b;c;z)\,,
\end{array} \right.
\ee
where $R(a,b)=-\Psi(a)-\Psi(b)-2 \gamma$.
The above asymptotic formula for the {\it zero-balanced} case $a+b=c$ 
is due to
Ramanujan (see  \cite{Be}). This formula is implied by 
\cite[15.3.10]{AS}. Note that $R(\frac{1}{2}, \frac{1}{2}) = \log 16$.

For complex $a,b,c$, and $z$, with $|z| < 1$, we let
\be \label{eq:Mdef}
\M(z) = \M(a,b,c,z) = 
z(1-z) \left( v_1(z) \frac{dv}{dz} - v(z)\frac{dv_1}{dz} \right) \,.
\ee
From (\ref{eq:hypder}) it is easy to see that
\bea 
\M & = & (c-a)(u v_1 + u_1 v) + (2(a-c)+b)v v_1 \label{eq:Mdef2} \\
 & = & (c-a)(u v_1 + u_1 v - v v_1) + (a+b-c)v v_1. \nonumber
\eea

It follows from
\cite[Corollary 3.13(5)]{AQVV} and (\ref{eq:gamref}) that

\be \label{eq:Maca}
\M(a,1-a,1,r) = \frac{1-a}{\Gamma(a)\Gamma(2-a)} =\frac{ \sin(\pi a) }{ \pi}
\ee
for $0 < a < 1$ and $0 \le r < 1$.
In particular,  we get the classical 
Legendre relation (\cite{AAR}, \cite{BF})
\be \label{eq:Legendre}
\M(1/2,1/2,1,r) = \frac{1}{\pi}\,.
\ee

The next result generalizes \cite[Theorem 3.9]{AQVV}.

\bigskip

\alku 
\theo \label{thm39aqvv} 
For $0 <a, b <c $, let the function $f$ be defined on $[0,\infty)$ by
$f(x) = F(a,b;c;1-e^{-x})$ and let $g(x) = f(x) \exp(-(a+b-c)x)$.
Then $f$ and $g$ are increasing, with $f(0)=g(0)=1$.  
If $a+b > c$, then $f(\infty) =\infty$. 
If $a+b = c$, then $f(\infty) =g(\infty) = \infty$. If $a+b < c$, then 
$f(\infty) = B(c,c-a-b)/B(c-a,c-b)$ and $g(\infty) = \infty$. 
Moreover, $h(x) = f'(x) e^{-x(a+b-c)}$ is also 
increasing, with $h(0) = ab/c$ and 
$h(\infty) = \Gamma(c) \Gamma(a+b+1-c)/(\Gamma(a) \Gamma(b))$
or $h(\infty) = \infty$ according as $a+b+1 > c$ or $a+b+1 \le c$. 
\loppu

\medskip

\proof 
The assertions
$f(0) = g(0) = 1$ and that $f$ is increasing, so that if $a+b < c$,
then $g$ is increasing, are all clear. In the three cases, 
$a+b < c, a+b = c$ and $a+b > c$,
the limiting values at $\infty$ are clear by (\ref{eq:hypas}).
Next, by (\ref{eq:hypas}),
$g(x) = F(c-a,c-b;c;1-e^{-x})$, which is clearly increasing.
Next, by differentiation,
$$
c/(ab) f'(x) = e^{-x} F(a+1,b+1;c+1; 1-e^{-x})\,.
$$
Hence, by (\ref{eq:hypas})
\begin{eqnarray*}
c/(ab) f'(x) & = & (e^{-x})(e^{-x(c-a-b-1)}) F(c-a,c-b;c+1; 1-e^{-x}) \\
& = & e^{x(a+b-c)} F(c-a,c-b;c+1; 1-e^{-x})\,,
\end{eqnarray*}
so that
$h(x) = (ab/c) F(c-a,c-b;c+1; 1-e^{-x})$, which is increasing, with 
boundary values $h(0) = ab/c$, and by 
(\ref{eq:hypas}) and (\ref{eq:gamapp}),
$$
h(\infty) =  (ab/c) 
\frac{ \Gamma(c+1) \Gamma(a+b+1-c)}{\Gamma(a+1)\Gamma(b+1)} 
= \frac{\Gamma(c)\Gamma(a+b+1-c)}{\Gamma(a)\Gamma(b)}\, , 
$$
if $a+b+1 > c$, and $= \infty$ if $a+b+1 \le c$. $\qquad \square$

\bigskip

\alku
\theo \label{th:aqvvg}
For $a, b, c, d > 0$, with $a+b > c > max \{a,b\}$, let the function $f$ 
be defined on $[0,\infty)$ by
$f(x) = F(a,b;c;1-(1+x)^{-1/d})$, and let 
$$
g(x) = (1+x)^{((c+d)-(a+b))/d} f'(x) \,.
$$ 
Then $g$ is increasing, with
$$
g(0)= \frac{ab}{cd} \:\:\:  \mbox{\rm and} \:\:\: 
g(\infty) = 
\frac{(a+b-c)\Gamma(c)\Gamma(a+b-c)}{d\, \Gamma(a)\Gamma(b)} \, .
$$
\loppu

\medskip

\proof With $u = 1-(1+x)^{-1/d}$, 
by (\ref{eq:hypder}) and (\ref{eq:hypas}), 
\begin{eqnarray*}
f'(x) & = &\frac{ab}{cd} (1+x)^{-1-(1/d)} F(a+1,b+1;c+1;u) \\
& = & \frac{ab}{cd} (1+x)^{(1/d)(a+b-c-d)} F(c-a,c-b;c+1;u)\,,
\end{eqnarray*}
so that
$g(x) = (ab/(cd)) F(c-a,c-b;c+1;u)$, which is clearly positive and 
increasing. The boundary value $g(0) = ab/(cd)$ is clear, while
the value of $g(\infty)$
follows from (\ref{eq:hypas}). $ \qquad \square$

\bigskip

\alku \rem \label{rem:triv} 
{\rm
(1) With $c = a+b$, Theorem \ref{thm39aqvv} reduces to 
\cite[Theorem 3.9]{AQVV}.

(2) For $a+b = c+d$, Theorem \ref{th:aqvvg} reduces to 
\cite[Theorem 3.10.]{AQVV}.
}
\loppu

\bigskip

\begin{figure}[htbp]
\begin{center}
\includegraphics[width=120mm]{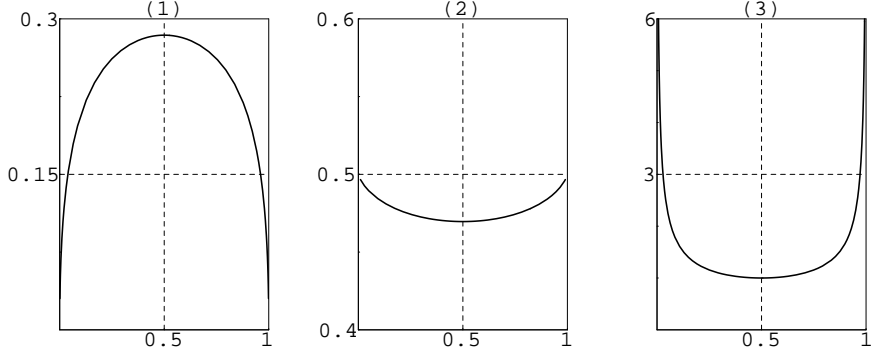}
\end{center}
\caption{(1) $\M(0.5, 1.0, 2.0, \cdot)$, 
(2) $\M(0.5, 1.0, 1.5, \cdot)$,
(3) $\M(0.5, 1.0, 1.0, \cdot)$
} \label{fig:mfunc}
\end{figure}

\bigskip

\alku \theo \label{Mtheo}
For positive constants $a,b,c$, the restriction to $(0,1)$ of the
continuous function $\M$ in (\ref{eq:Mdef}) has the 
following properties.
Denote $\M(x) = \M(a,b,c,x)$.

(1) $\M(x) = \M(1-x) > 0$ for all $x \in (0,1)$.

(2) If $a+b \le c$, then 
$\M(x)$ is bounded and extends continuously
to $[0,1]$. In particular, if $a+b=c=1$, then $\M(x)$ equals the 
constant $\sin(\pi a)/\pi$.

(3) If $a+b > c$, then $\M$ is unbounded on $(0,1)$, with
$\M(0+) = \M(1-) = \infty$.

(4) If $a+b < c < a+b+1,$ then 
$$
\lim_{ x \to 0+} x^{a+b-c} \M(x) = 
\lim_{x \to 1-} (1-x)^{a+b-c} \M(x) = 
\frac{\Gamma(c)\Gamma(a+b+1-c)}{\Gamma(a)\Gamma(b)}\,. 
$$

(5) If $a+b+1 = c$, then 
$$
\lim_{x \to 0+} \frac{\M(x)}{x \log(1/x)} = 
\lim_{x \to 1-} \frac{\M(x)}{(1-x) \log(1/(1-x))} = 
\frac{a+b}{B(a,b)}\,. 
$$

(6) If $a+b+1 < c$, then 
$$ 
\lim_{x \to 0+}\frac{\M(x)}{x} =
\lim_{x \to 1-}\frac{\M(x)}{1-x}=
\frac{ ab(2c-a-b-1)B(c,c-a-b)}{c(c-a-b-1)B(c-a,c-b)}\,. 
$$

(7) If $a+b > c$, then
$$ 
\lim_{x \to 0+} x^{a+b-c} \M(x) =
\lim_{x \to 1-} (1-x)^{a+b-c} \M(x) =
\frac{(a+b-c)B(c,a+b-c)}{B(a,b)}\,. 
$$
\loppu

\medskip

\proof (1) From (\ref{eq:fuv}) and (\ref{eq:Mdef}), we get
\begin{eqnarray*}
\M(x) & = & \frac{abx(1-x) }{c}(F(a+1,b+1;c+1;x)v(1-x) \\
      & & + F(a+1,b+1;c+1;1-x)v(x)) \\
 & = & G(x) + H(x) > 0\,. 
\end{eqnarray*}
Here  
$$
G(x) = \frac{abx(1-x)}{ c}F(a+1,b+1;c+1;x)v(1-x)
$$
and
$$
H(x) = \frac{abx(1-x) }{c} F(a+1,b+1;c+1;1-x)v(x) = G(1-x)\, .
$$

(2) First, if $a+b < c$, then from (\ref{eq:hypas}) and (\ref{eq:Mdef}),
$$
\M(0+) = \M(1-) = (c-a) u(1) + (a+b-c) v(1) = 0\,. 
$$
Next, if $a+b = c$, then from (\ref{eq:hypas}),
$G(0+) = 0$ and 
$$
H(x) = \frac{ab(1-x)}{c} F(a,b;a+b+1;1-x) v(x)\,,
$$ 
so that
$H(0+) = (ab/c) F(a,b;a+b+1;1)$. 
Next,
$H(1-) = 0$ and 
$$
G(x) = \frac{abx}{c} F(a,b;a+b+1;x) v(1-x)\,,
$$ 
so that $G(1-) = (ab/c) F(a,b;a+b+1;1)$. 
Thus,
$$ 
\M(0+) = \M(1-) = (ab/c) F(a,b;a+b+1;1) = 1/B(a,b) \, .
$$

(3) Let $a+b > c$. Then
$$
H(x) = \frac{ab}{c}(1-x)x^{c-a-b} F(c-a,c-b;c+1;1-x) v(x)\,,
$$ 
so that $H(0+) = \infty$. 
Similarly,
$G(1-) = \infty$.

(4) By (\ref{eq:hypas}), 
\begin{eqnarray*}
\frac{c}{ab} \M(x) & = & x(1-x)F(a+1,b+1;c+1;x)F(a,b;c;1-x;) \\
 & & + x^{c-a-b} (1-x) F(c-a,c-b;c+1;1-x)F(a,b;c;x)\,, 
\end{eqnarray*}
so that the result follows by (1).

(5) By (\ref{eq:hypas}), 
\begin{eqnarray*}
\frac{c \M(x)}{abx \log(1/x)} & = & x(1-x)F(a+1,b+1;c+1;x)F(a,b;c;1-x) \\ 
 & & + \frac{1-x}{\log(1/x)} F(a+1,b+1;c+1;1-x)F(a,b;c;x)\,,
\end{eqnarray*}
so that the result follows from (1) and (\ref{eq:hypas}).

(6) By (\ref{eq:hypas}),
\begin{eqnarray*}
\frac{c\M(x)}{abx} & = & (1-x)F(a+1,b+1;c+1;x)F(a,b;c;1-x) \\
 & & + (1-x)F(a+1,b+1;c+1;x)\,,
\end{eqnarray*}
so that by (1) and (\ref{eq:hypas}),
\begin{eqnarray*}
\lim_{x \to 0+} \frac{\M(x)}{x} & = & \lim_{ x \to 1-}{\M(x)/(1-x)} \\
 & = & \frac{ab}{c} \left( \frac{B(c,c-a-b)}{B(c-a,c-b)} + 
\frac{B(c+1,c-a-b-1)}{B(c-a,c-b)} \right) \\
 & = & \frac{ab(2c-a-b-1)}{(c(c-a-b-1)} \frac{B(c,c-a-b)}{B(c-a,c-b)}.
\end{eqnarray*}

(7) By (\ref{eq:hypas}), 
\begin{eqnarray*}
\frac{c}{ab} x^{a+b-c} \M(x) & = & x(1-x)F(a+1,b+1;c+1;x)F(a,b;c;1-x) \\
 & & + (1-x)F(c-a,c-b;c+1;1-x)F(a,b;c;x)\,,
\end{eqnarray*} 
so that by (1) and (\ref{eq:hypas}),
\begin{eqnarray*}
\lim_{x \to 0+} x^{a+b-c} \M(x) & = & 
\lim_{ x \to 1-} (1-x)^{a+b-c} \M(x) \\
 & = & \frac{ab}{c} F(c-a,c-b;c+1;1) \\
 & = & \frac{ab}{c} \frac{B(c+1,a+b+1-c)}{B(a+1,b+1)} \\
 & = & \frac{(a+b-c)B(c,a+b-c)}{B(a,b)} \,. \qquad \square 
\end{eqnarray*}

\bigskip

\alku \lem \label{Mdifflemma}
The function $\M$ in (\ref{eq:Mdef}) satisfies the differentiation formula
\begin{eqnarray*}
\frac{d \M(a,b,c,z)}{dz} & = & \frac{1}{z(1-z)}
\Bigl((c-a)[ (1-c+(a+b-1)z) u(z)v_1(z) \\ 
 & & + (-a-b+c+(a+b-1)z) u_1(z) v(z)] \\ 
 & & +(1-2z)[(c-a)(a+2b-1) -b^2]v(z) v_1(z)\Bigr)\,. 
\end{eqnarray*}
\loppu

\medskip

\proof 
Denote $D = d/dz$.
Then, by (\ref{eq:hypder}), (\ref{eq:hypder2}), and the chain rule, we get
$$
z Du = (a-1)(v-u), \:\: z(1-z) Dv = (c-a)u + (a-c+bz)v
$$
$$
(1-z) Du_1 = (1-a)(v_1 - u_1), \:\:
z(1-z) Dv_1 = - ((c-a)u_1 + (a-c+b(1-z))v_1).
$$
Hence, by the product rule, after simplification, we get
$$
z(1-z) D(uv_1)
= -(c-a)uu_1 + (a-1)(1-z)vv_1 + ((1-a)(1-z) +(c-a-b(1-z)))uv_1,
$$
$$
z(1-z) D(u_1 v)
= (c-a)uu_1 +(1-a)zvv_1 +(a-c +(a+b-1)z)u_1v,
$$
and
$$
z(1-z)D(vv_1)
= (c-a)(uv_1 - u_1v) + b(2z-1) vv_1. 
$$
Substituting these we get
\begin{eqnarray*}
z(1-z) D \M(a,b,c,z)
& = & (c-a)[(1-c+(a+b-1)z)u v_1 \\
& & + (-a-b+c +(a+b-1)z)u_1 v] \\
& & + (1-2z)[(c-a)(a+2b-1)-b^2] vv_1. \qquad \square
\end{eqnarray*}

\bigskip

\alku \rem \label{Mspec}
{\rm
If we put $c=a+b=1$ in Lemma \ref{Mdifflemma} above, then we get the familiar 
fact that
$$
\frac{d}{dz} \M(a,1-a,1,z) = 0.
$$
We also observe that $D \M(a,b,c,1/2) = 0$.
}
\loppu

\bigskip


\section{Generalized elliptic integrals} 
\label{sect:modular}

The following two important Theorems, \ref{monothos} and \ref{monotquot},
are indispensable in simplified proofs for monotonicity of the 
quotient of two functions. The first one, called
{\it L'H\^opital's Monotone Rule}, appears in 
\cite[Theorem 1.25]{AVV}, 
while a more general version of the second one appears in 
\cite{BK} and \cite[Lemma 2.1]{PV}.

\bigskip

\alku \theo \label{monothos}
Let $-\infty < a < b < \infty$, and let 
$f,\, g \, : \, [a,b] \to \R$ be continuous on $[a,b]$, 
differentiable on $(a,b)$. Let $g'(x) \neq 0$ on $(a,b)$. If
$f'(x)/g'(x)$ is increasing (decreasing) on $(a,b)$, then so are
$$
\frac{f(x)-f(a)}{g(x)-g(a)} \:\:\:\:\:  and \:\:\:\:\:
\frac{f(x)-f(b)}{g(x)-g(b)}.
$$
If $f'(x)/g'(x)$ is strictly monotone, then the monotoneity in the
conclusion is also strict.
\loppu

\bigskip

\alku \lem \label{monotquotapu}
Let $\{a_n\}$ and $\{b_n\}$ be real sequences with $b_n > 0$ for all $n$.
If the sequence $\{a_n/b_n\}$ is increasing (decreasing), then
$$
T_n = \sum_{k = 0}^{n} (n-k) (a_{n-k} b_k - a_k b_{n-k})
$$
is positive (negative) for $n = 1,2,\ldots$.
\loppu

\medskip

\proof
It is enough to prove the case $\{a_n/b_n\}$ is increasing, since the other
one is similar. Clearly
$$
T_1 = a_1 b_0 - a_0 b_1 = b_0 b_1 \left( \frac{a_1}{b_1} - 
\frac{a_0}{b_0} \right) > 0.
$$
Let $n \ge 2$. First let $n = 2m$ be even. Then
\begin{eqnarray*}
T_n = T_{2m} & = & 
 \sum_{k = 0}^{2m} (2m - k) (a_{2m-k} b_k - a_k b_{2m-k}) \\
 & = & \sum_{k = 0}^{2m-1} (2m - k)(a_{2m-k} b_k - a_k b_{2m-k}) \\
 & = & \sum_{k = 0}^{m-1} (2m - k)(a_{2m-k} b_k - a_k b_{2m-k}) + 0 \\
 & & + \sum_{k = m+1}^{2m-1} (2m - k)(a_{2m-k} b_k - a_k b_{2m-k}) \\
 & = & \sum_{k = 0}^{m-1} (2m - k)(a_{2m-k} b_k - a_k b_{2m-k}) \\
 & & + \sum_{k = 1}^{m-1} k (a_k b_{2m-k} - a_{2m-k} b_k) \\
 & = & 2 m (a_{2m} b_0 - a_0 b_{2m}) + 
 \sum_{k = 1}^{m-1} (2m - 2k)(a_{2m-k} b_k - a_k b_{2m-k}) \\
 & = & 2 m b_0 b_{2m} \left(\frac{a_{2m}}{b_{2m}} - \frac{a_0}{b_0} \right)
 + \sum_{k=1}^{m-1}(2m - 2k) b_k b_{2m-k} 
 \left(\frac{a_{2m-k}}{b_{2m-k}} - \frac{a_k}{b_k}\right) > 0.
\end{eqnarray*}
Next, let $n = 2m + 1$ be odd. Then
\begin{eqnarray*}
T_n & = & 
 \sum_{k = 0}^{2m+1} (2m+1-k) (a_{2m+1-k} b_k - a_k b_{2m+1-k}) \\
 & = & \sum_{k = 0}^{2m} (2m+1-k)(a_{2m+1-k} b_k - a_k b_{2m+1-k}) \\
 & = & (2m+1)(a_{2m+1} b_0 - a_0 b_{2m+1}) +  
 \sum_{k = 1}^{m} (2m+1-k)(a_{2m+1-k} b_k - a_k b_{2m+1-k}) \\
 & & + \sum_{k = m+1}^{2m} (2m+1-k)(a_{2m+1-k} b_k - a_k b_{2m+1-k}) \\
 & = & (2m+1)(a_{2m+1} b_0 - a_0 b_{2m+1}) 
 + \sum_{k = 1}^{m} (2m+1-k)(a_{2m+1-k} b_k - a_k b_{2m+1-k}) \\
 & & + \sum_{k = 1}^{m} k (a_k b_{2m+1-k} - a_{2m+1-k} b_k) \\
 & = & (2m+1)(a_{2m+1} b_0 - a_0 b_{2m+1}) +  
 \sum_{k = 1}^{m} (2m+1-2k)(a_{2m+1-k} b_k - a_k b_{2m+1-k}) \\
 & = & (2m+1)(b_0 b_{2m+1}) \left(\frac{a_{2m+1}}{b_{2m+1}} - 
 \frac{a_0}{b_0} \right) \\
 & & + \sum_{k=1}^{m}(2m+1-2k) b_k b_{2m+1-k} 
 \left(\frac{a_{2m+1-k}}{b_{2m+1-k}} - \frac{a_k}{b_k}\right) > 0.
\qquad \square
\end{eqnarray*}

\bigskip

\alku \theo \label{monotquot}
Let $\sum_{n = 0}^{\infty} a_n x^n$ and $\sum_{n = 0}^{\infty} b_n x^n$
be two real power series converging on the interval $(-R, R)$. If 
the sequence $\{a_n/b_n\}$ is increasing (decreasing), and $b_n > 0$ for
all $n$, then the function
$$
f(x) = \frac{\sum_{n=0}^{\infty} a_n x^n}{\sum_{n=0}^{\infty} b_n x^n}
$$
is also increasing (decreasing) on $(0,R)$. In fact, the function
$f'(x) \left(\sum_{n=0}^{\infty} b_n x^n \right)^2$ has positive
Maclaurin coefficients.
\loppu

\medskip

\proof
\begin{eqnarray*}
f'(x) \left(\sum_{n=0}^{\infty} b_n x^n \right)^2 
 & = & \sum_{n = 0}^{\infty} b_n x^n \sum_{n=0}^{\infty} n a_n x^{n-1} -
\sum_{n = 0}^{\infty} a_n x^n \sum_{n=0}^{\infty} n b_n x^{n-1} \\
 & = & (1/x) \sum_{n=0}^{\infty} \left( \sum_{k = 0}^{n} 
(n-k)(a_{n-k} b_k - a_k b_{n-k}) \right) x^n.
\end{eqnarray*}
The result follows from Lemma \ref{monotquotapu}. $\qquad \square$

\bigskip

The following theorem solves the corresponding problem in
the case where we have a quotient of two polynomials instead of
two power series.

\bigskip

\alku \theo \label{monpolquot}
Let $f_n(x) = \sum_{k=0}^{n} a_k x^k$ and
$g_n(x) = \sum_{k=0}^{n} b_k x^k$ be real polynomials, with
$b_k > 0$ for all $k$. If the sequence $\{a_k/b_k\}$ is 
increasing (decreasing), then so is the function
$f_n(x)/g_n(x)$ for all $x > 0$. In fact, 
$g_n f_n' - f_n g_n'$ has positive (negative) coefficients.
\loppu

\medskip

\proof
We prove the increasing case by induction on $n$. The proof of the 
decreasing case is similar. Let first $n = 1$. Then
$$
\frac{f_1(x)}{g_1(x)} = \frac{a_0 + a_1 x}{b_0 + b_1 x}.
$$
Hence
\begin{eqnarray*}
g_1 f_1' - f_1 g_1' & = & (b_0 + b_1 x)a_1 - (a_0 + a_1 x)b_1 \\
 & = & a_1 b_0 - a_0 b_1 \\
 & = & b_0 b_1 \left( \frac{a_1}{b_1} - \frac{a_0}{b_0} \right) > 0.
\end{eqnarray*}

Next, assume that the claim holds for all $k \le n$. Now
$$
\frac{f_{n+1}(x)}{g_{n+1}(x)} = 
\frac{f_n(x) + a_{n+1} x^{n+1}}{g_n(x) + b_{n+1} x^{n+1}}.
$$
We get
\begin{eqnarray*}
g_{n+1} f_{n+1}' - f_{n+1} g_{n+1}' 
& = &  (g_n + b_{n+1} x^{n+1})(f_n' + (n+1)a_{n+1}x^n) \\
 & & - (f_n + a_{n+1} x^{n+1})(g_n' + (n+1)b_{n+1}x^n) \\
& = & (g_n f_n' - f_n g_n') + 
      (n+1) x^n (g_n a_{n+1} - f_n b_{n+1}) \\
 & & + x^{n+1} (b_{n+1} f_n' - a_{n+1} g_n') \\
& = & (g_n f_n' - f_n g_n') + 
      (n+1) x^n \sum_{k=0}^{n}(a_{n+1} b_k - b_{n+1} a_k) x^k \\
 & & + x^{n+1} \sum_{k=1}^{n} k (a_k b_{n+1} - b_k a_{n+1}) x^{k-1} \\
& = & (g_n f_n' - f_n g_n') + 
      \sum_{k=0}^{n}(n+1-k)(a_{n+1}b_k - b_{n+1}a_k) x^{n+k} \\
& = & (g_n f_n' - f_n g_n') + 
      \sum_{k=0}^{n}(n+1-k) b_k b_{n+1}
      \left(\frac{a_{n+1}}{b_{n+1}} - \frac{a_k}{b_k} \right) x^{n+k}. \\
\end{eqnarray*}
Hence each coefficient is positive. $\qquad \square$

\bigskip

\alku \lem \label{muabcmonot}
Let $a,b,c, K > 0$.

(1) If $a+b \ge c$, then $\mu_{a,b,c} : (0,1) \to (0, \infty)$ and
$\varphi_K^{a,b,c} : (0,1) \to (0,1)$ are decreasing and
increasing homeomorphisms, respectively.

(2)  If $a+b > c$, then the function 
$f(r) = (r/r')^{2(a+b-c)} \mu_{a,b,c}(r)$
is decreasing from $(0,1)$ onto $((B(a,b)^2)/(2B(c,a+b-c)), B(c,a+b-c)/2)$.

(3) If $a+b > c$ and $K > 1$, then
$$
r < \varphi_K^{a,b,c}(r) < K^{1/(2(a+b-c))} r \,,
$$ 
for all $r \in (0,1)$.

(4) If $a + b < c$, then 
$\mu_{a,b,c} : [0,1] \to [B(a,b)/(2d), B(a,b) d/2]$,
is a decreasing homeomorphism, where $d= F(a,b;c;1)$ is given by 
(\ref{eq:hypas}). In this case, $\varphi_K^{a,b,c}$ is defined if and
only if $K \ge 1$. If $K > 1$, then $\varphi_K^{a,b,c}$ maps 
$[0,1]$ onto a proper subset of $[0,1]$.
\loppu

\medskip

\proof
(1) Since $a,b,c > 0$ and $a+b \ge c$, it follows 
from (\ref{eq:hypas}) that the function
$F(a,b;c;r)$ is an increasing homeomorphism of
$[0,1)$ onto $[1,\infty)$. Hence
the function $\mu_{a,b,c}$ is a strictly decreasing 
homeomorphism $(0,1) \to (0, \infty)$. The assertion about
$\varphi_K^{a,b,c}$ follows from these facts.

(2) From (\ref{eq:hypas}), 
$$ 
f(r) = \frac{B(a,b)}{2} 
\frac{F(c-a,c-b;c;{r'}^2)}{F(c-a,c-b;c;r^2)}, ,
$$ 
which is decreasing with required limiting values at $0$ and $1$.

(3) With $s = \varphi_K^{a,b,c}(r) > r$, from (2) we get
\begin{eqnarray*}
f(s) & = & (s/s')^{2(a+b-c)} \mu_{a,b,c}(s) =  
(s/s')^{2(a+b-c)} \mu_{a,b,c}(r)/K \\
 & < & f(r) = (r/r')^{2(a+b-c)} \mu_{a,b,c}(r)\, ,
\end{eqnarray*}
so that
$$
(s/s')^{2(a+b-c)} < K (r/r')^{2(a+b-c)}\,,
$$
which gives 
$s/s' < K^{1/(2(a+b-c))} r/r'$. Hence
$$
r < s < K^{1/(2(a+b-c))} r(s'/r') < K^{1/(2(a+b-c))} r\,.
$$ 

(4) For $a,b,c > 0$ and $a+b < c$, (\ref{eq:hypas}) implies that
$F(a,b;c;r)$ is an increasing homeomorphism 
$[0,1] \to [1, F(a,b;c;1)]$.
Now the claim follows from (\ref{eq:muacdef}) and (\ref{eq:phiackdef}).
$\qquad \square$

\bigskip

\begin{figure}[htb]
\begin{center}
\includegraphics[width=70mm]{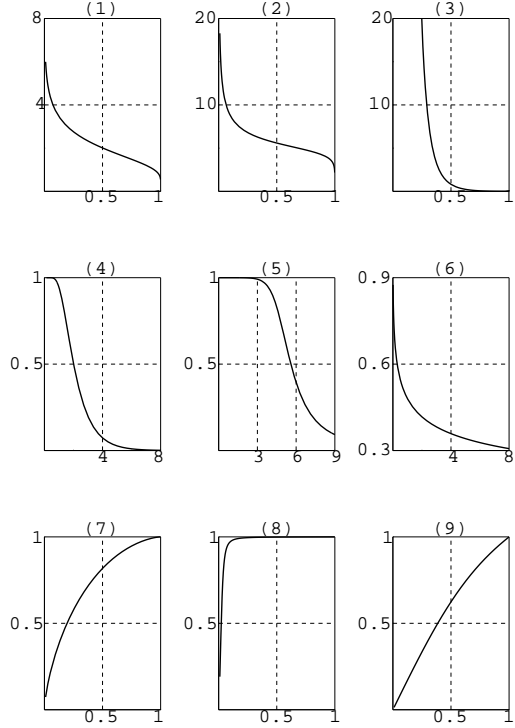}
\end{center}
\caption{(1) $\mu_{0.5,0.5,1}$, 
(2) $\mu_{0.1,1.1,1}$,  
(3) $\mu_{2.5,1.5,2}$, 
(4) $\mu^{-1}_{0.5,0.5,1}$,
(5) $\mu^{-1}_{0.1,1.1,1}$, \newline 
(6) $\mu^{-1}_{2.5,1.5,2}$,
(7) $\varphi_{1.5}^{0.5,0.5,1}$,
(8) $\varphi_{2.5}^{0.1,1.1,1}$, 
(9) $\varphi_{3.5}^{2.5,1.5,2}$.
} \label{fig:funcs1}
\end{figure}

\bigskip

The function $\mu_{a,b,c}$ is a natural generalization for the
function $\mu_a$ in \cite[(1.3)]{AQVV}. Namely,
\be \label{eq:mua1mua}
\mu_{a,1-a,1}(r) = \frac{B(a,1-a)}{2} 
\frac{F(a,1-a;1;{r'}^2)}{F(a,1-a;1;r^2)} =
\mu_a(r)\,,
\ee
since by (\ref{eq:gamref}) and (\ref{eq:betadef}),
$$
B(a,1-a) = \frac{\Gamma(a) \Gamma(1-a)}{\Gamma(1)} = 
\frac{\pi}{\sin(\pi a)}\,.
$$

Clearly (\ref{eq:muacdef}) and (\ref{eq:phiackdef}) imply the identities 
(cf. \cite[(4.8), (4.9)]{AQVV}),
\be \label{eq:muprod} 
\mu_{a,b,c}(r) \mu_{a,b,c}(r') = \left(\frac{B(a,b)}{2}\right)^2\,,
\:\: r \in (0,1)\,, 
\ee
\be \label{eq:musum}
\mu_{a,b,c}^{-1}(x)^2 + \mu_{a,b,c}^{-1}(y)^2 = 1\,,
\ee
where $x,y > 0$ with $xy = (B(a,b)/2)^2$, and
\be \label{eq:phiacksum} 
\varphi^{a,b,c}_K(r)^2 + \varphi^{a,b,c}_{1/K}(r')^2 = 1\,.
\ee

Moreover, from (\ref{eq:muacdef}) and (\ref{eq:hypas}) we get,
for $c < a + b$, 
\be \label{eq:mutrans} 
{r'}^{2(c-a-b)} \mu_{a,b,c}(r)=
r^{2(c-a-b)} \mu_{c-a,c-b,c}(r)\,, \:\: r \in (0,1)\,.
\ee

\bigskip

\alku \lem \label{distflem}
Let $f$ be a bijection from a real interval $I$ onto
$(0,\infty)$ and let $g$ be defined on $I$ by $g(x) = a f(x)$, where 
$a > 0$ is a
constant. Then for each constant $K > 0$, we have
$$
f^{-1}(K f(x)) = g^{-1}(K g(x))\,.
$$
\loppu

\medskip

\proof Let $u = f^{-1}(K f(x))$. Then $f(u) = K f(x)$, so that
$a f(u) = a K f(x)$, that is, $g(u) = K g(x)$. Hence 
$u = g^{-1}(K g(x))$. $\qquad \square$

\bigskip

\alku \rem \label{rem:phitil}
{\rm
For $a,b,c > 0$ with $a+b \ge c$, $K > 0$, 
and $r \in (0,1)$ denote
$$
\mut_{a,b,c}(r) = \frac{F(a,b;c;{r'}^2)}{F(a,b;c;r^2)}
$$
and
$$
\phit^{a,b,c}_{K}(r) = \mut_{a,b,c}^{-1}(\mut_{a,b,c}(r)/K)\,.
$$
By Lemma \ref{distflem} we see that
$$
\phit^{a,b,c}_{K}(r) = \varphi^{a,b,c}_{K}(r)\,.
$$
}
\loppu

\bigskip

\alku \rem \label{rem:eklimits} {\rm Observe first that
$ \lim_{c \to a+} B(a,c-a) = \infty$
and $\Gamma(c-a) (c-a,n) = \Gamma(c-a+n)$, which
tends to $\Gamma(n) = (n-1)!$, as $c \to a+$.

(1) $\lim_{c \to a+} \left({ \K}_{a,c}(r) - (B(a,c-a)/2) \right)
 = \log (1/r')$.

(2) $\lim_{ c \to a+} \left({ \E}_{a,c}(r) - (B(a,c-a)/2) \right)
= \log (1/r')  -
 \left( \frac{1}{2} \sum_{n= 1}^{\infty} { \frac{ r^{2n}} {a+n-1} } 
\right)$.

In particular, for each fixed $r \in (0,1)$
all the three functions $\K_{a,c}(r), \E_{a,c}(r)$, and 
$\E_{a,c}(r) - {r'}^2 \K_{a,c}(r)$
tend to $\infty$ as $c \to a+$.
}
\loppu

\bigskip

\alku \theo \label{EKdifftheo}
The following differentiation formulae hold:
\be \label{eq:Kdiff}
\frac{d \K_{a,b,c}}{dr} = 
\frac{2}{r {r'}^2} \left((c-a)\E_{a,b,c}+(b {r}^2+a-c) \K_{a,b,c} \right)\,,
\ee
\be \label{eq:Ediff}
\frac{d \E_{a,b,c}}{dr} =
\frac{2(a-1)}{r}\left(\K_{a,b,c}-\E_{a,b,c} \right) \,,
\ee
\bea 
\frac{d}{dr}(\K_{a,b,c}-\E_{a,b,c}) & = &  
\frac{2}{r {r'}^2} \left(((c-a)-(1-a){r'}^2)\E_{a,b,c} \right. 
\label{eq:KmEdiff} \\
&& \left. +((a+b) {r}^2-c +{r'}^2) \K_{a,b,c} \right)\,, \nonumber
\eea 
\be \label{eq:KmrpEdiff}
\frac{d}{dr} (\E_{a,b,c} - {r'}^2 \K_{a,b,c}) = 
\frac{2}{r}((1-c)\E_{a,b,c} + (c-1-(b-1)r^2)\K_{a,b,c})\,,
\ee
\be \label{eq:muacdiff}
\frac{d}{dr} \mu_{a,b,c}(r) = 
- \frac{B(a,b)\M(a,b,c,r^2)}{r{r'}^2 v(r^2)^2} =
- \frac{B(a,b)^3\M(r^2)}{4r{r'}^2 \K^2}\,,
\ee
\be \label{eq:phiacdiff}
\frac{\M(a,b,c,s^2)}{\M(a,b,c,r^2)} \frac{ds}{dr} = 
\frac{1}{K} 
\frac{s {s'}^2 v(s^2)^2}{r {r'}^2 v(r^2)^2} = 
\frac{1}{K} 
\frac{s {s'}^2 \K(s)^2}{r {r'}^2 \K(r)^2}\,, 
\:\: s = \varphi^{a,b,c}_{K}(r)\,.
\ee
\loppu

\medskip

\proof
From (\ref{eq:hypder})
$$ \frac{dF}{dz}= \frac{(c-a)u +(a-c+bz)v}{z(1-z)} \,.
$$
Put $z=r^2$ and multiply both sides by $B(a,b)/2\,.$ 
This gives (\ref{eq:Kdiff}). 

For (\ref{eq:Ediff}) recall that by (\ref{eq:hypder2})
$$ \frac{du}{dz}= \frac{(a-1)}{z}(v-u) \,. $$
Put $z=r^2$ and multiply both sides by $B(a,b)/2\,.$
This gives (\ref{eq:Ediff}).

The formulae (\ref{eq:KmEdiff}) and (\ref{eq:KmrpEdiff}) 
follow from (\ref{eq:Kdiff}) and
(\ref{eq:Ediff}).
  
From (\ref{eq:Mdef}), putting $z = r^2$, we get
$$
\M(r^2) = -r^2 {r'}^2 \frac{d \mu}{dr} 
\left(\frac{2}{B(a,b)}\right)^3 \K^2 \left(\frac{1}{2r}\right),
$$
so that
$$
\frac{d \mu}{dr} = - \frac{B(a,b)^3 \M(r^2)}{4 r {r'}^2 \K^2}.
$$

Denote $s = \varphi^{a,b,c}_{K}(r)$. Then by (\ref{eq:phiackdef})
$$
\mu_{a,b,c}(s) = \frac{1}{K} \mu_{a,b,c}(r)
$$
so that
$$
\frac{d}{ds} \mu_{a,b,c}(s) \frac{ds}{dr} = 
\frac{1}{K} \frac{d}{dr} \mu_{a,b,c}(r)\,. 
$$
Now (\ref{eq:phiacdiff}) follows from (\ref{eq:muacdiff}). 
$\qquad \square$

\bigskip

\alku \lem \label{EKmonot} (cf.\ \cite[Lemmas 5.2 and 5.4]{AQVV})
For $0< a, b < \min \{c,1\}$ and $c \le a+b \, ,$  denote $\K = \K_{a,b,c}$
and $\E = \E_{a,b,c}$. Then the function

(1) $f_1(r) = (\K - \E)/(r^2 \K)$ is strictly 
increasing from $(0,1)$ onto $(b/c, 1)$. In
particular, we have the sharp inequality,
$$
\frac{b}{c} < \frac{\K - \E}{r^2 \K} < 1
$$
for all $r \in (0,1)$.

(2) $f_2(r) = (\E-{r'}^2 \K)/r^2$ has positive Maclaurin
coefficients and maps $(0,1)$ onto 
\newline
$(B(a,b)(c-b)/(2c), C)$, where
$$ 
 C = \frac{B(a,b)B(c,c+1-a-b)}{2B(c+1-a,c-b)} \, .
$$

(3) $f_3(r) = {r'}^{-2(c+1-a-b)} \E$ 
has positive Maclaurin coefficients and maps $[0,1)$ onto 
\newline
$(B(a,b)/2, \infty)$.  

(4) $f_4(r) = {r'}^{2(a+b-c)} \K$ 
has positive Maclaurin coefficients and maps
$[0,1)$ onto \newline 
$[B(a,b)/2, B(c,a+b-c)/2)$.

(5) $f_5(r) = {r'}^{-2} \E$
has positive Maclaurin coefficients and maps
$[0,1)$ onto 
\newline
$[B(a,b)/2, \infty)$.  

(6) $f_6(r) = {r'}^2 \K$ 
has negative Maclaurin coefficients, except for the constant term, and maps
$[0,1)$ onto $(0,B(a,b)/2]$.

(7) $f_7(r) = \K$ has positive Maclaurin coefficients and is 
log-convex from $[0,1)$ onto 
\newline
$[B(a,b)/2, \infty)$.  In fact, $(d/dr)(\log \K) $ also has
positive Maclaurin coefficients. 

(8) $f_8(r) = (\E - {r'}^2 \K)/(r^2 \K)$ 
is strictly decreasing from $(0,1)$ onto $(0,1-(b/c))$.

(9) $f_9(r) = (\K - \E)/(\E - {r'}^2 \K)$ is
strictly increasing from $(0,1)$ onto $(b/(c-b), \infty)$.

(10) $f_{10}(r)= (\K - B(a,b)/2)/ \log(1/r')$ 
is strictly increasing from
$(0,1)$ onto 
\newline
$(ab B(a,b)/c, D)$, where $D=1$ if $c=a+b$ and 
$D= \infty$ if $c <a+b$.

(11) $f_{11}(r) = (B(a,b)/2 - {r'}^2 \K)/r^2$
has positive Maclaurin coefficients and maps
$(0,1)$ onto 
$(B(a,b)(c-ab)/(2c), B(a,b)/2)$.

(12) $f_{12}(r) = (\K-B(a,b)/2)/({r'}^{2(c-a-b)}-1)$,  
for $a+b>c$, is strictly increasing 
from $(0,1)$ onto $(abB(a,b)/(2c(a+b-c)), B(c,a+b-c)/2)$.

(13) $f_{13}(r) = ((1-a-(b-c){r}^2)\E - 
(1-a){r'}^2 \K)/r^2$
has negative Maclaurin coefficients, except for the constant term,
with
$$  
f_{13}(0+)=(c+1-a)(c-b)B(a,b)/(2c)
$$
and $f_{13}(1-) = (c+1-a-b) \E(1)$.

(14) $f_{14}(r) =  (c-a)\E - (b-a){r'}^2\K$
has negative Maclaurin coefficients, except for the constant term,
and maps $[0,1]$ onto $[C,D]$, where
$$
 C = (c-a) \E(1) \: \:
\mbox{\rm and} \:\:
D = (c-b)B(a,b)/2.
$$
\loppu

\medskip

\proof
(1) From (\ref{eq:hypder})
$$
\frac{dF(a^-)}{dz} =\frac{(a-1)}{z} \left( F -F(a^-)\right)= 
\frac{(a-1)b}{c} F(a,b+1;c+1;z) \,.
$$ 
Putting $z= r^2$ and multiplying by $B(a,b)/2$, we get
$$ 
\frac{\K -\E}{r^2} = \frac{b B(a,b)}{2c} F(a,b+1;c+1;r^2) \,. 
$$
Hence
$$
f_1(r) = \frac{b}{c} \frac{F(a,b+1;c+1;r^2)}{F(a,b;c;r^2)} \, .
$$
Thus $f_1(0)=b/c$. The ratio of the coefficients of the numerator and
denominator equals
$$
\frac{b}{c} \frac{(a,n) (b+1,n) (c,n)}{(c+1,n)(a,n)(b,n)}= 
\frac{b}{c}\frac{b+n}{b}\frac{c}{c+n}
=\frac{b+n}{c+n}= 1- \frac{c-b}{c+n},
$$
which is increasing in $n$. Hence the result follows from Theorem 
\ref{monotquot}.
The limit $f_1(1-)= 1$ follows from
(\ref{eq:hypas}).

(2) From  (\ref{eq:hypder}) and (\ref{eq:hypder2}), we get
$$  
\frac{dF(a-)}{dz} =  \frac{(a-1)}{z} (F- F(a-))= 
(a-1)(F +  \frac{(b-c)}{c} F(c+)) \,.
$$
Hence putting $z= r^2$, we get
$$
\frac{\K- \E}{r^2}= \K +  \frac{(b-c) B(a,b)}{2c} F(a,b;c+1;r^2) \,.
$$
Thus, 
$$
f_2(r)=  \frac{(c-b)B(a,b)}{2c} F(a,b;c+1;r^2)  \, ,
$$
which proves the assertion. The limiting values follow from (\ref{eq:hypas}).

(3) By (\ref{eq:hypas}),
\begin{eqnarray*}
f_3(r) & = & \frac{B(a,b)}{2} \frac{F(a-1,b;c;r^2)}{{r'}^{2(c+1-a-b)}} \\
 & = & \frac{B(a,b)}{2} F(c+1-a, c-b; c; r^2),
\end{eqnarray*}
and the result follows.

(4) By (\ref{eq:hypas}),
\begin{eqnarray*}
f_4(r) & = & \frac{B(a,b)}{2} {r'}^{2(a+b-c)} F(a,b;c;r^2) \\
 & = & \frac{B(a,b)}{2} F(c-a, c-b; c; r^2),
\end{eqnarray*}
and the result follows.

(5) By (\ref{eq:hypas})
$$ 
f_5(r)=  \frac{B(a,b)}{2 {r'}^2} F(a-1,b;c;r^2)=
\frac{B(a,b)}{2} (1-r^2)^{-(a+b-c)} F(c+1-a,c-b;c;r^2)\, ,
$$
which is a product of two Maclaurin series with positive 
coefficients, hence the result.

(6) We have that 
$$
f_6(r) = \E(r) - (\E(r) - {r'}^2 \K(r))\,. 
$$
Hence, by (2), $f_6(r) -B(a,b)/2$ has all Maclaurin coefficients negative.

(7) The positivity of Maclaurin 
coefficients and the limiting values are clear. 
Next, by (\ref{eq:Kdiff}), after simplification, we get
$$
\frac{d}{dr} \log \K
= 
\frac{2r}{{r'}^2 \K} \left( (a+b-c) \K +  
\frac{ (c-a)(\E - {r'}^2 \K)}{r^2} \right) \, ,
$$
which  has positive Maclaurin coefficients by (2) and (6).

(8) $f_8(r) = 1 - f_1(r)$, so that the result follows from (1).

(9) $f_9(r) = f_1(r)/f_8(r)$, hence the result
follows from (1) and (8). 

(10) The ratio of the coefficients equals
$$ 
\frac{B(a,b) (a,n)(b,n)}{(c,n)(n-1)!} \equiv T_n \, ,
$$
say. Then
$$
\frac{T_{n+1}}{T_n} = \frac{(a+n)(b+n)}{(c+n)n} >1 \, ,
$$
so that the monotonicity follows from \ref{monotquot}. 
The limiting values follow from (\ref{eq:hypas}).

(11) 
\begin{eqnarray*}
f_{11}(r) & = & \K - \frac{\K -B/2}{r^2} \\
 & = & \frac{B}{2} \left(1 -\frac{ab}{c} + 
\sum_{n=1}^{\infty} \left(\frac{(a,n)(b,n)}{(c,n)n !} 
(1- \frac{(a+n)(b+n)}{(c+n)(n+1)}) r^{2n}\right) \right),
\end{eqnarray*}
which has all coefficients positive. The limiting values are obvious.

(12) The ratio has the indeterminate form $0/0$ at $r = 0$. The derivative
ratio equals
$$
\frac{ab B(a,b)}{2 c (a+b-c)} 
\frac{F(a+1,b+1;c+1;r^2)}{{r'}^{2(c-a-b-1)}} =
\frac{ab B(a,b)}{2 c (a+b-c)} F(c-a,c-b;c+1;r^2),
$$
by (\ref{eq:hypas}), and so the result follows by Lemma
\ref{monothos}.

(13) From (\ref{eq:hypder})
and (\ref{eq:hypder2})
\begin{eqnarray*}
\frac{dF(a-)}{dz} & = & 
\frac{1}{1-z}\left((a-1+b-c)F(a-) + \frac{(c+1-a)(c-b)}{c}F(a-,c+)\right) \\ 
 & = & \frac{(a-1)}{z}(F-F(a-)) \, .
\end{eqnarray*}
Multiplying by $z(1-z)B(a,b)/2$, we get
$$
z(a-1+b-c) \E+ \frac{(c+1-a)(c-b)}{c}z F(a-,c+) \frac{B(a,b)}{2} =
(a-1)(1-z)(\K-\E) \, .
$$
With $z=r^2$ this yields the result.

(14) In the Gauss contiguous relation in \cite[Exercise 21(9), p.\ 71]{R1},  
if we change $z$ to $r^2$ and multiply by $B(a,b)/2$, 
then we get $f_{14}(r) =  ((c-b)B(a,b)/2) F(a,b-1;c;r^2)$,
so that the assertions on the coefficients follow.
The limiting values are clear.  $\qquad \square$

\bigskip

\alku \rem \label{rem:417-11}
{\rm 
In the classical case, $a=b=1/2$ and $c = 1$, the boundary 
values in the above result \ref{EKmonot} (13) reduce to
$f_{13}(0+) = 3\pi/8$ and $f_{13}(1-) = 1$, 
showing that the above result is quite sharp. 
}
\loppu

\bigskip

\alku \theo \label{avvtheo1_52}
Let $0 < a < c \le 1$, $b = c-a$,   
$R = R(a,c-a) = -\Psi(a) - \Psi(c-a) - 2 \gamma$, and 
$B = B(a,c-a)$. Then

(1) The function $f(r) = \K_{a,c}(r) + \log r'$ 
has negative Maclaurin coefficients, except for the constant term, 
and maps $[0,1)$ onto $(R/2, B/2]$.

(2) The function $g(r) = \K_{a,c}(r) + (1/r^2) \log r'$ 
has positive Maclaurin coefficients and maps
$(0,1)$ onto $((B-1)/2, R/2)$
if $a, b \in (0,1)$, 
while it has negative Maclaurin coefficients and maps
$(0,1)$ onto $(R/2, (B-1)/2)$ if $a, b \in (1,\infty)$.

(3) The function $h(r) = r^2 \K_{a,c}(r)/\log(1/r')$ 
is strictly decreasing (respectively, increasing)  from $(0,1)$ 
onto $(1,B)$ if $a,b \in (0,1)$  (respectively, onto $(B,1)$, if 
$a, b \in (1,\infty) $). 

(4)  The function $k(r) =  \K_{a,c}(r) / \log ((e^{R/2})/ r')$ is 
strictly decreasing from $(0,1)$ onto $(1, B/R)$. 
\loppu

\medskip
 
\proof 
(1) That $f(0+) = B(a,c-a)/2$ is clear and
$f(1-) = R(a,c-a)/2$ follows from \cite[Theorem 1.52 (2)]{AVV}.
Next,
$$
2 f(r) = B + B \sum_{n=1}^{\infty}\left(\frac{(a,n) (b,n)}{(c,n) n!}
- \frac{1}{n}\right) r^{2n}.
$$
Thus, we need to show that $T_n = (a,n) (b,n)/((c,n) (n-1)!) < 1$.
Now
$$
\frac{T_{n+1}}{T_n} = \frac{(a+n) (b+n)}{(c+n) n}
= \frac{(a+b)n + n^2 + ab}{(a+b)n + n^2} > 1,
$$
and $\lim_{n \to \infty} T_n = 1$ by Stirling's formula. Hence
$T_n < 1$ for all $n = 1,2,3,\ldots$.

(2)
$$
g(0+) = \frac{B}{2} + \lim_{r \to 0} \frac{\log r'}{r^2} 
= \frac{B}{2} - \frac{1}{2}.
$$
Next, $g(r) = f(r) + (1/{r}^2 - 1) \log r'$, so that
$g(1-) = f(1-) = R/2$, from (1). Next,
\begin{eqnarray*}
2 g(r) & = & B F(a,b;c;r^2) - \frac{1}{r^2} \log \frac{1}{1-r^2} \\
 & = & \sum_{n= 0}^{\infty} \frac{1}{n+1} 
\left( \frac{B (a,n) (b,n) (n+1)}{(c,n) n!} -1 \right) r^{2n}.
\end{eqnarray*}
Let 
$$
T_n = \frac{B (a,n) (b,n) (n+1)}{(c,n) n!}.
$$
Then
$$
\frac{T_{n+1}}{T_n} = \frac{(a+n)(b+n)(n+2)}{(c+n)(n+1)^2}.
$$
Now $(a+n)(b+n)(n+2) - (c+n)(n+1)^2 = -(1-ab)n - (a+b-2ab)$, which
is negative (positive) if $a,b \in (0,1)$ ($a,b \in (1, \infty)$).
By Stirling's formula, $\lim_{n \to \infty} T_n = 1$. Hence the result follows.

(3) 
$$
h(r) = \frac{g(r)}{\frac{1}{r^2} \log \frac{1}{r'}} + 1, 
$$
so that $h(0+) = B$ and $h(1-) = 1$ are clear. Next,
$$
h(r) = \frac{B \sum_{n=0}^{\infty} \frac{(a,n)(b,n)}{(c,n) n!} r^{2n}}{
\sum_{n=0}^{\infty}\frac{1}{n+1}r^{2n}},
$$
so that the coefficient ratio equals $(a,n) (b,n) (n+1)/((c,n) n!)$,
which is decreasing if $a,b \in (0,1)$ and increasing if $a,b \in (1,\infty)$.
Hence the result follows from \cite[Lemma 2.1]{PV}.

(4)  We have $k(r) = 1 + (f(r) - R/2)/ \log(e^{R/2} /r')$. Hence the 
result follows from (1). $\qquad \square$

\bigskip

\alku \theo \label{kmudiff}
(1) Let $0 < a, \, b < c$ and $2ab < c \le a+b < c + 1/2$.
Then the function $f(r) = r' \K(r)$ 
is strictly decreasing from $[0,1)$ onto $(0, B(a,b)/2]$.

(2) Let $0 < a, b < c < a+b$. Then the function
$g(r) = {r'}^{2(a+b-c)}(\K - \E)/r^2$ has positive Maclaurin
coefficients and maps $(0,1)$ onto $(bB(a,b)/(2c), B(c,a+b-c)/2)$.

(3) Let $0 < a, \, b < 1$, $c = a+b$, and $a(2b+1) < b+1 < 1/a$. 
Then the function
$h(r) = (\K - \E)/\log(1/r')$ is decreasing from $(0,1)$ onto
$(1, bB(a,b)/c)$.
\loppu

\medskip

\proof
(1) 
Clearly $f(0) = B(a,b)/2$ and by (\ref{eq:hypas}), $f(1) = 0$.
We have $f(r) = g(r)/h(r)$, where
$$
g(r) = B(a,b) \sum_{0}^{\infty}\frac{(a,n)(b,n)}{(c,n) n!} r^{2n}\,, 
$$
and
$$
h(r) = 2 \sum_{0}^{\infty} d_n r^{2n}\,, 
$$
with
$d_0 = 1$ and 
$d_n = (1 \cdot 3 \cdot \dots \cdot (2n-1))/
(2 \cdot 4 \cdot \dots \cdot (2n))$
for $n = 1,2,3,\ldots$.
Hence, the coefficient ratio equals
$$
T_n = \frac{B(a,b)(a,n)(b,n)2^n}{(c,n) 
\cdot 1 \cdot 3 \cdot \dots \cdot(2n-1)}\,. 
$$
Then
$$
\frac{T_{n+1}}{T_n} = \frac{2(n+1)(n+b)}{(2n+1)(n+c)}\,. 
$$
Now
$$
(2n+1)(n+c)- 2(n+a)(n+b) = n(2c+1-2a-2b) + c-2ab  > 0\,, 
$$
so that $T_n$ is decreasing and hence by Theorem \ref{monotquot},
$f$ is also decreasing.

(2) From (\ref{eq:hypder}), 
(\ref{eq:hypas}) and (\ref{eq:Ediff}), we get
$$ 
f(r) = \frac{bB(a,b)}{2c} {r'}^{2(a+b-c)} F(a,b+1;c+1;r^2)
= \frac{bB(a,b)}{2c} F(c+1-a,c-b;c+1;r^2). 
$$
Hence the assertion follows from (\ref{eq:hypas}).

(3) As in (2), from (\ref{eq:hypder}), 
(\ref{eq:hypas}) and (\ref{eq:Ediff}), we get
$$
h(r) = \frac{bB(a,b)}{c} \frac{r^2 F(a,b+1;c+1;r^2)}{2\log(1/r')}. 
$$
Writing the Maclaurin
series expansion of both the numerator and the denominator, the ratio of
coefficients equals
$$
T_n = \frac{bB(a,b)}{c} \frac{(a,n) (b+1,n) (n+1)}{(c+1,n) n!}
= \frac{B(a,b)(a,n)(b,n+1)(n+1)}{(c,n+1) n!}. 
$$
Hence
$$
\frac{T_{n+1}}{T_n} = \frac{(a+n)(b+n+1)(n+2)}{(c+n+1)(n+1)^2}. 
$$
Then,
$$
(c+n+1)(n+1)^{2} - (a+n)(b+n+1)(n+2) = n(1-a-ab) + (1+b-2ab-a) > 0.
$$
Hence, $T_n$ is
decreasing, so that the result follows by (\ref{eq:hypas}) and 
Theorem \ref{monotquot}. 
$\qquad \square$

\bigskip

\alku \rem \label{v040602rmk}
{\rm 
Theorem \ref{kmudiff}(3) generalizes \cite[Lemma 5.2 (12)]{AQVV}. The 
latter follows if we put $c = 1$ in Theorem \ref{kmudiff}(3). 
}
\loppu

\bigskip

\ots{Differential equations}
The hypergeometric function 
$w = F(a,b;c;z)$ 
satisfies the differential equation \cite[p.\ 54]{R1}
$$
z(1-z)w'' + (c-(a+b+1)z)w' - abw = 0\,.
$$
Changing the variable $z$ to $z^2$, this reduces to
$$
z(1-z^2) w'' + (2c-1-(2a+2b+1)z^2)w' - 4abzw = 0\,.
$$
In particular, the generalized elliptic integrals,
$w = \K_{a,c}(r)$ and $w = \K_{a,c}'(r)$, satisfy, respectively,
the differential equations
\be \label{eq:Kdy}
r {r'}^2 w'' + (2c - 1 - (2c+1)r^2) w' - 4a(c-a) r w = 0\,,
\ee
and
\be \label{eq:Kpdy}
r {r'}^2 w'' - (1 - (2c+1)r^2)w' - 4a(c-a) r w = 0\,.
\ee
In the special case $c = 1$ the above two equations coincide 
(cf.\ \cite[(4.3)]{AQVV}, \cite[3.8.19, p.\ 75]{L}).
Next, the generalized elliptic integrals
$w = \E_{a,c}(r)$ and $w = \E_{a,c}'(r)$ satisfy, respectively,
the differential equations
\be \label{eq:Edy}
r {r'}^2 w'' - (2c-1){r'}^2 w' + 4(1-a)(c-a) r w = 0\,,
\ee
and
\be \label{eq:Epdy}
r {r'}^2 w'' - (1+(2c-1)r^2) w' + 4 (1-a)(c-a) r w = 0\,.
\ee
In the special case $c = 1$ the equations (\ref{eq:Edy}) and
(\ref{eq:Epdy}) are still different unlike in the case of 
(\ref{eq:Kdy}) and (\ref{eq:Kpdy}) 
(cf.\ \cite[3.8.17, p.\ 74 and 3.8.23, p.\ 75]{L}).

\bigskip

\ots{Correction} In \cite[(4.3), p.\ 14]{AQVV}
the first differential equation
has a symmetry property, namely, it is satisfied both by $\K_{a}$
and $\K'_{a}$.
However, the second differential equation is satisfied only by 
$\E_{a}$, and {\em not} by $\E'_{a}$.
The differential equation satisfied by $w = \E'_{a}$ is obtained by
putting $c = 1$ in (\ref{eq:Epdy}). Thus it is
$$
r{r'}^2 w'' - (1+r^2) w' + 4(1-a)^2 r w = 0\,.
$$

\bigskip

We use the notation
$$ 
S_{w} = \left( \frac{w''}{w'} \right)^{'}  - \frac{1}{2} 
\left( \frac{w''}{w'} \right)^2
$$
for the {\it Schwarzian derivative}.

\bigskip

\alku \lem \label{Rainlem} \cite[p.\ 9]{R2} 
Let $w_1, w_2$ be linearly independent
solutions of the differential equation $w''+p(z) w' + q(z) w = 0$.
Then $W= w_1/w_2$ satisfies the  differential equation 
$$ 
S_W(z) = 2 q(z) -p'(z) - p(z)^2/2 \, .
$$
\loppu

\bigskip

\alku \theo \label{gmudeqn} Let $0<a,b<1$ and $2c=a+b+1$. Then the modulus
$\mu = \mu_{a,b,c}$ satisfies the differential equation
$$ 
S_{\mu}(r) = 2 q(r) -p'(r) - p(r)^2/2\,,
$$
where
$$ 
p(r)=  \frac{2c-1-(4c-1) r^ 2}{ r {r'}^ 2} \, , 
\quad q(r) = -\frac{4ab}{{r'}^2 } \, .
$$
\loppu

\medskip

\proof Follows immediately from Lemma \ref{Rainlem} and (\ref{eq:Kdy}).
$\qquad \square$

\bigskip

\alku \theo \label{geqn} Let $0<a,b<1$ and $2c=a+b$. Then the function
$\nu = \E'/\E$ satisfies the differential equation
$$  
S_{\nu}(r) = 2 q(r) -p'(r) - p(r)^2/2 \,,
$$
where
$$ 
p(r)=  \frac{2c-1-(4c-1) r^ 2}{ r {r'}^ 2} \, , 
\quad q(r) = -\frac{4(a-1)b}{{r'}^2 } \, .
$$
\loppu

\medskip

\proof Follows immediately from Lemma \ref{Rainlem} and (\ref{eq:Edy}).
$\qquad \square$

\bigskip











\bigskip

\alku \lem \label{lem:430} 
If $0 < a < \min \{c,1\}$ and $c \le a+(1/2),$ then the 
function $f(r) = r {\K}_{a,c}(r) /{\rm arth} r $
is strictly decreasing 
from $(0,1)$ onto $( 1,B/2),$ where $B = B(a,c-a).$
\loppu

\medskip

\proof
Let $f(r) = g(r)/h(r)$, where
$g(r) = r \K_{a,c}(r)$ and $h(r) = \mbox{\rm arth}\, r$.
Then $g(0) = h(0) = 0$ and
$$
g'(r)/h'(r) = 2(c-a) \E_{a,c}(r) + (1-2(c-a))({r'}^2) \K_{a,c}(r)\,,
$$
which, by Lemma \ref{EKmonot} (6), being a sum of two
strictly decreasing functions, is also so. Hence, the monotonicity 
follows from L'H\^opital's Monotone Rule, Lemma \ref{monothos}.
Finally, $f(0+) = B/2$ by L'H\^opital's Rule and 
$f(1-) = g'(1-)/h'(1-) = 1$, again 
by L'H\^opital's Rule. $\qquad \square$

\bigskip

{\bf Acknowledgments.} The authors wish to thank the Finnish National
Academy of Science and Letters, the Finnish Mathematical Society, and the
Departments of Mathematics at the University of Helsinki and the
University of Turku for generous support. The authors are indebted to
Heikki Ruskeep\"a\"a for providing expert help with the 
$\mbox{\rm Mathematica}^{\circledR}$ software. They are also grateful
to D.\ Askey, B.\ C.\ Carlson, N.\ Stylianopoulos, and T.\ Sugawa for
interesting discussions and correspondence.

\bigskip

{\small
\noindent
{\sc V.\ Heikkala}\\
SSH Communications Security Corp.\\
Valimotie 17\\
FIN--00380 Helsinki\\
Finland\\
email: {\tt ville.heikkala@ssh.com} \\
fax: +358-20-5007051 \\
\\
\\
{\sc M.K.\ Vamanamurthy } \\
Department of Mathematics \\
The University of Auckland \\
P B 92019, Auckland \\
New Zealand \\
email: {\tt vamanamu@math.auckland.ac.nz} \\
fax: +64-9-3737457 \\
\\
\\
{\sc M.\ Vuorinen}\\
Department of Mathematics \\
FIN-20014 University of Turku \\
FINLAND \\
email: {\tt vuorinen@utu.fi} \\
fax: +358-2-3336595
}

\bigskip

{\small
\begin{thebibliography}{AQVV}

\bibitem [AS] {AS} {\sc M.\ Abramowitz and I.\ A.\ Stegun, editors}:
{\em Handbook of Mathematical Functions with Formulas, Graphs and
Mathematical Tables}, Dover, New York, 1965.

\bibitem [Ah] {Ah} {\sc L.\ V.\ Ahlfors}:
{\em Complex Analysis}, 2nd ed., McGraw-Hill, New York, 1966.

\bibitem [AQVV] {AQVV} {\sc G.\ D.\ Anderson, S.\ -L.\ Qiu, 
M.\ K.\ Vamanamurthy, and M.\ Vuorinen}:
{\em Generalized elliptic integrals and modular equations},
Pacific J.\ Math.\ 192 (2000), 1--37.

\bibitem [AVV] {AVV} {\sc G.\ D.\ Anderson, M.\ K.\ Vamanamurthy, and
M.\ Vuorinen}: {\em Conformal Invariants, Inequalities, and Quasiconformal
Maps}, J.\ Wiley, 1997.

\bibitem [AAR] {AAR} {\sc G.\ E.\ Andrews, R.\ Askey, and
R.\ Roy}: {\em Special Functions},
Cambridge Univ.\ Press, 1999.

\bibitem[B1]{B1} \textsc{A.  Baricz:} Landen-type inequality 
for Bessel function, \emph{  Comput. Methods Funct. Theory} 
\textbf{5}  (2005),  no. 2, 373--379. 


\bibitem[B2]{B2} \textsc{A. Baricz:} Functional inequalities 
involving special functions,
 \emph{ J. Math. Anal. Appl.} \textbf{319} (2006), 450--459.

\bibitem[B3]{B3} \textsc{A. Baricz:} Functional inequalities 
involving special functions II,  
\emph{ J. Math. Anal. Appl.} \textbf{327} (2007), 1202--1213.

\bibitem[B4]{B4} \textsc{A. Baricz:} Tur\'an type inequalities 
for generalized complete elliptic integrals,
\emph{ Math. Z.} 256 (2007), 895--911.




\bibitem [Be] {Be} {\sc B.\ C.\ Berndt}:
{\em Modular equations in Ramanujan's lost notebook},
Number theory, 55--74, Trends Math., Birkh\"auser, Basel, 2000.
                                                                              
\bibitem [BBG] {BBG} {\sc B.\ C.\ Berndt, S.\ Bhargava and F.\ G.\ Garvan}:
{\em Ramanujan's theories of elliptic functions to alternative bases}, 
Trans.\ Amer.\ Math.\ Soc., 347 (1995), 4163--4244.


\bibitem [BK] {BK} {\sc M.\ Biernacki and J.\ Krzy\.{z}}: 
{\em On the monotonicity of 
certain functionals in the theory of analytic functions}, 
Ann.\ Univ.\ M.\ Curie-Sklodowska, 2 (1955), 134--145.

\bibitem [BB] {BB} {\sc J. M. Borwein and P. B. Borwein}: {\em Pi and
the AGM}, John Wiley \& Sons, New York, 1987.

\bibitem [Bo] {Bo} {\sc F.\ Bowman}: {\em Introduction to Elliptic 
Functions with applications,} English Universities Press Ltd., London, 
1953.

\bibitem [BF] {BF} {\sc P.\ F.\ Byrd and M.\ D.\ Friedman}:
{\em Handbook of Elliptic Integrals for Engineers and 
Scientists}, 2nd ed., Grundlehren Math.\ Wiss. Vol.\ 67,
Springer-Verlag, Berlin, 1971.

\bibitem [CLT] {CLT} {\sc H.\ H.\ Chan, W-C.\ Liaw and V.\ Tan}: 
{\em Ramanujan's class invariant $\lambda_n$ and a new class of 
series for $1/\pi$},  
J.\ London Math.\ Soc.\ (2) 64 (2001), no.\ 1, 93--106.

\bibitem [DV] {DV} {\sc V. N. Dubinin and M. Vuorinen}: 
{\em On conformal moduli of polygonal quadrilaterals}, 
 arXiv:math.CV/0701387

\bibitem [HLVV] {HLVV} {\sc
 V. Heikkala, H. Lind\'en, M. K. Vamanamurthy, and  M. Vuorinen}:
Generalized elliptic integrals II,  arXiv:math.CA/0701438


\bibitem [Hen] {Hen} {\sc P.\ Henrici}: 
{\em Applied and Computational Complex Analysis},
Vol.\ III,  Wiley, New York, 1986.

\bibitem [Her] {Her} {\sc J.\ Hersch}:
{\em On harmonic measures, conformal moduli and some elementary symmetry
methods},  
J.\ Analyse Math.\ 42 (1982/83), 211--228.

\bibitem [KS] {KS} {\sc D. Karp and S. M. Sitnik}:
Inequalities and monotonicity of ratios for generalized
hypergeometric function, arXiv:math.CA/0703084.

\bibitem [L] {L} {\sc D.\ F.\ Lawden}:
{\em Elliptic Functions and Applications}, 
Applied Math.\ Sciences, Vol.\ 80, Springer-Verlag, New York, 1989.

\bibitem [LV] {LV} {\sc O. Lehto and K. I. Virtanen}: 
{\em Quasiconformal Mappings in the Plane,} 2nd ed., Grundlehren Math. 
Wiss., Band 126, Springer-Verlag, New York, 1973.

\bibitem [M] {M} {\sc A.\ I.\ Markushevich}:
{\em Theory of Functions of a Complex Variable},
Vol.\ II, Prentice-Hall, Englewood Cliffs, NJ, 1965.


\bibitem [N] {N} {\sc Z.\ Nehari}:
{\em Conformal mapping},
McGraw-Hill Book Co., Inc., New York, Toronto, London, 1952. viii+396 pp.

\bibitem [PV] {PV} {\sc S.\ Ponnusamy and M.\ Vuorinen}:
{\em Asymptotic expansions and inequalities for hypergeometric 
functions}, Mathematika 44 (1997), 278--301.

\bibitem [Q] {Q} {\sc S.\ L.\ Qiu}:
{\em Gr\"otzsch ring and Ramanujan's modular equations},
(Chinese) Acta Math.\ Sinica 43 (2000), no.\ 2, 283--290.

\bibitem [QV1] {QV1} {\sc S.\ L.\ Qiu and M.\ Vuorinen}:
{\em Duplication inequalities for the ratios of hypergeometric functions}, 
Forum Math.\ 12 (2000), no.\ 1, 109--133.

\bibitem [QV2] {QV2} {\sc S.-L. Qiu and M. Vuorinen}:
{\em  Landen inequalities for hypergeometric functions,} 
Nagoya Math. J. 154 (1999), 31-56.

\bibitem [R1] {R1} {\sc E.\ D.\ Rainville}:
{\em Special Functions}, Macmillan, New York, 1960.

\bibitem [R2] {R2} {\sc E.\ D.\ Rainville}:
{\em Intermediate differential equations}, 2nd ed.  Macmillan, New York, 1964.

\bibitem [RV] {RV} {\sc A. Rasila and M. Vuorinen}: 
{\em Experiments with the moduli of quadrilaterals}, 
arXiv:math.NA/0703149,
 Rev.\ Roumaine Math.\ Pures Appl.\ 51 (2006), 747--757,  



\bibitem [S] {S} {\sc L-C.\ Shen}:
{\em On an identity of Ramanujan based on the hypergeometric series
${}_2 F_1(1/3, 2/3;1/2;x)$},
J.\ Number Theory 69 (1998), no.\ 2, 125--134.


\bibitem[WZQC]{WZQC} \textsc{G. Wang, X. Zhang, S.-L. Qiu, and Y. Chu:} 
The bounds of the solutions to generalized modular equations,  
\emph{J. Math. Anal. Appl.} \textbf{321} (2006), 589--594.


\bibitem [WW] {WW} {\sc E.\ T.\ Whittaker and G.\ N.\ Watson}:
{\em A Course of Modern Analysis}, 4th ed., Cambridge Univ.\ Press, London,
1927.



\bibitem[ZWC]{ZWC} \textsc{X. Zhang, G. Wang, and Y. Chu:} 
Some inequalities for the
generalized Gr\"otzsch function,  \emph{ Proc. Edinburgh
Math. Soc.} ( to appear )



\end{thebibliography}
}
\end{document}